\documentclass{amsart}
\usepackage{amsfonts}
\usepackage{amssymb}
\input epsf 
\newtheorem{thm}{Theorem}[section]
\newtheorem{cor}[thm]{Corollary}

\newtheorem{prop}[thm]{Proposition}
\theoremstyle{definition}
\newtheorem{define}[thm]{Definition}
\theoremstyle{remark}

\numberwithin{equation}{section}

\begin{document}

\centerline{Appeared in: Compositio Mathematica, 101 (1996), 179-215}

\vspace{1cm}

\title[Arnold-Liouville with singularities]{Symplectic topology of integrable
Hamiltonian systems, I: Arnold-Liouville with
singularities}

\author{Nguyen Tien Zung}
\address{SISSA - ISAS, Via Beirut 2-4, Trieste 34013, Italy (old address)}
\email{tienzung@math.univ-montp2.fr {\it URL}: www.math.univ-montp2.fr/\~{}tienzung}
\date{Ref. S.I.S.S.A. 153/94/M, Oct. 1994. Revised Jan 1995, accepted in final form 15 Feb. 1995.
This version (2001) contains a list of errata}

\begin{abstract}
The classical Arnold-Liouville theorem describes the geometry of an integrable
Hamiltonian system near a regular level set of the moment map. Our results describe
it near a nondegenerate singular level set: a tubular neighborhood of a connected
singular nondegenerate level set, after a normal finite covering, admits a
non-complete system of action-angle functions (the number of action functions is
equal to the rank of the moment map), and it can be decomposed topologically,
together with the associated singular Lagrangian foliation, to a direct product of
simplest (codimension 1 and codimension 2) singularities. These results are
essential for the global topological study of integrable Hamiltonian systems.
\end{abstract}

\keywords{Integrable Hamiltonian system, nondegenerate singularity, torus action,
topological decomposition, action-angle coordinates.}

\subjclass{58F07, 58F14, 58F05, 70H05}

\maketitle

\section{Introduction}

The classical Arnold-Liouville theorem \cite{Arnold} describes the geometry of an
integrable Hamiltonian system (IHS) with $n$ degrees of freedom near regular level
sets of the moment map: locally near every regular invariant torus there is a system
of action-angle coordinates, and the corresponding foliation by invariant Lagrangian
tori is trivial. In particular, there is a free symplectic $n$-dimensional torus
action which preserves the moment map. The local base space of this trivial
foliation is a disk equipped with a unique natural integral affine structure, given
by the action coordinates. A good account on action-angle coordinates, and their
generation for Mishchenko-Fomenko noncommutative integrable systems and Libermann
symplectically complete foliations, can be found in the paper \cite{DD} by Dazord
and Delzant, and references therein.

However, IHS's met in classical mechanics and physics always
have singularities, and most often these singularities are nondegenerate in
a natural sense. Thus the study of nondegenerate singularities
(=tubular neighborhoods of singular level sets of the moment map) is one of
the main step toward understanding the topological structure of general
IHS's.

The question is: What can be said about the structure of nondegenerate
singularities of IHS's?
What remains true from Arnold-Liouville theorem in this case?

In recent years, several significant results have been obtained concerning this
question. Firstly, the local structure of nondegenerate singularities is now
well-understood, due to the works by Birkhoff, Williamson, R\"ussmann, Vey and
Eliasson (cf. e.g. \cite{Birkhoff,Williamson,Russmann,Vey,Eliasson}). The second
well-understood thing is elliptic singularities and systems ``de type torique'',
through the works by Dufour, Molino and others (cf. e.g.
\cite{BM,Dufour,DM,Eliasson,Molino,Veeravalli}). In particular, Boucetta and Molino
\cite{BM} observed that the orbit space of the integrable system in case of elliptic
singularities is a singular integral affine manifold with a simple singular
structure. Integrable systems with elliptic singularities are closely related to
torus actions on symplectic manifolds, which were studied by Duistermaat, Heckman
and others (cf. e.g. \cite{DH,Audin,Delzant}). Third, codimension 1 nondegenerate
singularities were studied by Fomenko \cite{Fomenko} from the topological point of
view, and was shown to admit a system of $(n-1)$ action and $(n-1)$ angle
coordinates in my thesis \cite{ZungThesis}. Fourth, some simple cases of codimension
2 singularities were studied by Lerman, Umanskii and others (cf. e.g.
\cite{LU1,LU2,LU3,Bolsinov,Zou}). However, for general nondegenerate singularities,
the question had remained open. Here we listed only ``abstract'' works. There are
many authors who studied singularities of concrete well-known integrable systems,
from different viewpoints, see e.g.
\cite{AS,BolsinovP,CK,Devaney,Gavrilov,Fomenko2,GOC,Kharlamov,Oshemkov,ZungLada}.

The aim of this paper is to give a more or less satisfactory answer
to the question posed above, which incorporates cited results.
Namely, we will show that,
for a nondegenerate singularity of codimension $k$
in an IHS with $n$ degrees of freedom there is still a
locally free symplectic torus group action (of dimension $n-k$),
which preserves the moment map. There is a
normal finite covering of the singularity, for which the above torus action
becomes free, and the obtained torus
bundle structure is topologically trivial.
There is a coisotropic section to this trivial bundle, which gives rise to a
set of $(n-k)$ action and $(n-k)$ angle coordinates. And the most important
result, which completes the picture is as follows: the singular Lagrangian
foliation associated to a nondegenerate singularity is, up to a normal
finite covering, homeomorphic to a direct product of singular Lagrangian
foliations of simplest (codimension 1 and/or focus-focus codimension 2)
singularities.

Thus, we have the following topological classification of nondegenerate
singularities of integrable Hamiltonian systems: every such singularity has a unique
natural\footnote{erratum: The uniqueness is true when $k=n$, i.e. when the
singularity contains a fixed point, but when $k < n$ in some cases it's false. When
it's false, the number of minimal models is finite} {\sc canonical model} (or may be
also called minimal model because of its minimality), which is a direct-product
singularity with a free component-wise action of a finite group on it, cf. Theorem
\ref{thm:adp}. This finite group is a topological invariant of the singularity, and
may be called Galois group.

Canonical models for geodesic flows on multi-dimensional ellipsoids were
computed in \cite{ZungSphere}. In \cite{BauZung} we will discuss how the
method of separation of variables allows one to compute canonical model
of singularities of integrable systems. It is an interesting problem to
study singularities of systems which are integrable by other methods
(Lie-group theoretic, isospectral deformation, etc.).
In this direction, some results connected with Lie algebras were
obtained by Bolsinov \cite{BolsinovP}. Isospectral
deformation was deployed by Audin for some systems
\cite{AudinQ,AS,AudinTop}.

All the results listed above are true under a mild condition of topological
stability. This condition is satisfied for all nondegenerate singularities
of all known IHS's met in mechanics and physics. (Notice that this condition
is rather a structural stability but not a stability in Lyapunov's sense).
Without this condition, some of the above results still hold.

We suspect that our results (with some modifications) are also valid for integrable
PDE. In this direction, some results were obtained by Ercolani and McLaughlin
\cite{EM}. (Note: There seems to be an inaccuracy in Theorem 1 of \cite{EM}).

The above results have direct consequences in the global topological study of IHS's,
which are discussed in \cite{ZungSurgery}. For example, one obtains that the
singular orbit space of a nondegenerate IHS is a stratified integral affine manifold
in a natural sense, and one can compute homologies by looking at the bifurcation
diagram, etc. Our results may be useful also for perturbation theory of nearly
integrable systems (see e.g. \cite{Koiller} for some results concerning the use of
the topology of singularities of IHS's in Poincar\'e-Melnikov method), and for
geometric quantization.

{\bf About the notations}. In this paper, ${\bf T}$ denotes a torus, $\bf D$
a disk, ${\bf Z}$
the integral numbers, $\bf R$ the real numbers, $\bf C$ the complex numbers,
${\bf Z}_2$ the two element group. All IHS's are assumed to be smooth,
nonresonant, and by abuse of language, they are identified with Poisson
actions generated by $n$ commuting first integrals.  We refer to
\cite{DD,Nekhoroshev}
for a treatment of Arnold-Liouville theorem, and when we say that we use
this theorem we will often in fact use the following statement: If $X$ is a
vector field with periodic flow such that on each Liouville torus it equals
a linear combination of Hamiltonian vector fields of first integrals,
then $X$ is a symplectic vector field.

{\bf Organization of the paper}. In Section \ref{section:local}
we recall the main known results about the local structure of nondegenerate
singularities of IHS's. Section \ref{section:singleaf}
is devoted to some a priori
non-local properties of nondegenerate singularities and singular Lagrangian
foliations associated to IHS's. Then in
Section \ref{section:codim1} and Section \ref{section:focus} we describe the
structure of nondegenerate singularities in two simplest cases: codimension
1 hyperbolic and codimension 2 focus-focus. The main results, however, are
contained in the subsequent Sections, where we prove the existence of
Hamiltonian torus actions in the general case,
topological decomposition, and action-angle
coordinates. In the last section we briefly discuss various notions of
nondegenerate IHS's, and the use of integrable surgery in obtaining
interesting symplectic structures and IHS's.

{\bf Acknowledgements}. This paper is a detailed and corrected version of the
announcement \cite{Zung94}, and is based in part on my thesis \cite{ZungThesis}
defended at IRMA, Strasbourg, May/1994. I am indebted to Mich\`ele Audin, Anatoly T.
Fomenko, Alberto Verjovsky and the jury of my thesis for their support and
encouragement. I would like to thank also Alexey Bolsinov, Nicole Desolneux-Moulis,
Jean-Paul Dufour, Lubomir Gavrilov, Boris Kruglikov, Pierre Molino, and Jean-Claude
Sikorav, for useful discussions about the subject.

\section{Local structure of nondegenerate singularities}
\label{section:local}
Throughout this paper, the word IHS, which stands for integrable Hamiltonian
system, will always mean a $C^{\infty}$ smooth
Poisson ${\bf R}^n$ action in a smooth symplectic
$2n$-dimensional manifold $(M^{2n}, \omega)$, generated by a moment map ${\bf
F}: M \to {\bf R}^n$. We will always assume that the connected components of
the preimages of the moment map are compact.

Denote by $Q(2n)$ the set of quadratic forms on ${\bf R}^{2n}$. Under the
standard Poisson bracket, $Q(2n)$ becomes a Lie algebra, isomorphic to the
symplectic algebra $sp(2n, {\bf R})$. Its Cartan subalgebras have dimension
$n$.

Assume there is given an IHS with the moment map
${\bf F} = (F_1,...,F_n): M^{2n} \to {\bf R}^n$. Let $X_i = X_{F_i}$ denote
the Hamiltonian vector fields of the components of the moment map.

A point $x \in M^{2n}$ is called a {\sc singular point} for the above IHS,
if it is singular for the moment map: the differential
$D{\bf F} = (dF_1,...,dF_n)$ at $x$ has rank less than $n$. The {\sc corank}
of this singular point is ${\rm corank} D{\bf F} = n - {\rm rank} D{\bf F},
D{\bf F} = (dF_1,...,dF_n)$.

Let $x$ be a singular point as above. Let ${\mathcal K}$ be the kernel of
$D{\bf F}(x)$ and let ${\mathcal I}$
be the space generated by $X_i(x) (1 \leq i \leq n)$.
${\mathcal I}$ is a maximal isotropic subspace of ${\mathcal K}$ with
respect to the symplectic structure $\omega_0 = \omega (x)$. Hence the
quotient space ${\mathcal K} / {\mathcal I} = \overline{R}$ carries a natural
symplectic structure $\overline{\omega_0}$. $\overline{R}$ is symplectically
isomorphic to a subspace $R$ of $T_{x} M$ of dimension $2k$, $k$ being the
corank of the singular point. The quadratic parts of $F_1,...,F_n$ at
$x$ generate a subspace ${\mathcal F}^{(2)}_{R}(x)$ of the space of quadratic
forms on $R$. This subspace is a commutative subalgebra under the Poisson
bracket, and is often called in the literature the
{\sc transversal linearization} of ${\bf F}$.
With the above notations, the following definition is standard
(see e.g. \cite{Desolneux}).

\begin{define}
\label{define:singularpoint}
A singular point $x$ of corank $k$ is called {\sc nondegenerate} if
${\mathcal F}^{(2)}_{R}(x)$ is a Cartan subalgebra of the algebra of quadratic
forms on $R$. 
\end{define}

The following classical result of Williamson \cite{Williamson}
classifies linearized
nondegenerate singular points (see also \cite{Arnold}).

\begin{thm}[Williamson]
\label{thm:Williamson}
For any Cartan subalgebra ${\mathcal C}$ of Q(2n) there is a symplectic system
of coordinates $(x_1,...,x_n,y_1,...y_n)$ in ${\bf R}^{2n}$ and a basis
$f_1,...,f_n$ of ${\mathcal C}$ such that each $f_i$ is one of the following

$ \; f_i = x_i^2 + y_i^2$  (elliptic type)

$ \; f_i = x_iy_i$  (hyperbolic type)

$ \left.   \begin{array}{ll} f_i & = x_i y_{i+1}- x_{i+1} y_i \\
                            f_{i+1} & = x_i y_i + x_{i+1} y_{i+1}
           \end{array} \right\}$ (focus-focus type)
\end{thm}

For example, when $n=2$ there are 4 possible combinations: 1) $f_1, f_2$
elliptic; 2) $f_1$ elliptic, $f_2$ hyperbolic; 3) $f_1, f_2$ hyperbolic and
4) $f_1, f_2$ focus-focus. In some papers, these 4 cases are also called
center-center, center-saddle, saddle-saddle and focus-focus respectively
(cf. \cite{LU1}).

The above theorem is nothing but the classification up to conjugacy of
Cartan subalgebras of real symplectic algebras. In his paper, Williamson
also considered other subalgebras. Recall that a complex symplectic
algebra has only one conjugacy class of Cartan subalgebras. That is to say,
if we consider systems with analytic coefficients and complexify them, then
the above 3 types become the same.

\begin{define}
If the transversal linearization of an IHS at a nondegenerate
singular point $x$ of corank $k$ has
$k_e = k_e (x)$ elliptic components, $k_h = k_h (x)$ hyperbolic components,
$k_f = k_f (x)$ focus-focus
components ($k_e,k_h,k_f \geq 0, k_e + k_h + 2 k_f = k$),
then we will say that point $x$ has {\sc Williamson type} $(k_e,k_h,k_f)$.
Point $x$ is called {\sc simple} if $k_e + k_h + k_f =1$, i.e. if $x$ is
either codimension 1 elliptic, codimension 1 hyperbolic, or codimension 2
focus-focus.
\end{define}

It is useful to know the group of local linear automorphisms for simple
fixed points of linear IHS's (with one or two degrees of freedom
respectively), for these local linear automorphisms serve as components of
the linearization of local automorphisms of nonlinear IHS's. We have:

\begin{prop}
\label{prop:localautogroup}
Let $x$ be a simple fixed point of a linear IHS ($k_e + k_h + k_f =  1$). \\
a) If $x$ is elliptic then the group of local linear symplectomorphisms
which preserve the moment map is isomorphic to a circle ${\bf S}^1$. \\
b) If $x$ is hyperbolic then this group is isomorphic to ${\bf Z}_2 \times
{\bf R}^1$ = the multiplicative group of nonzero real numbers.  \\
c) If $x$ is focus-focus then this group is isomorphic to ${\bf S}^1 \times
{\bf R}^1$ = the multiplicative group of nonzero complex numbers.
\end{prop}

{\it Proof}.
We use the canonical coordinates given in Williamson's Theorem
\ref{thm:Williamson}. In elliptic case the group of local linear
automorphisms is just the group of rotations of the plane $\{x_1,y_1\}$.
In hyperbolic case, the group of endomorphisms of the plane $\{x_1, y_1\}$
which preserve the function $x_1y_1$ consists of elements
$(x_1,y_1) \mapsto (ax_1, y_1/a)$ and
$(x_1,y_1) \mapsto (ay_1, x_1/a)$ ($a \in {\bf R} \setminus \{0\}$).
But the symplectic condition implies that only the maps
$(x_1,y_1) \mapsto (ax_1, y_1/a)$ belong to our local automorphism group.
In focus-focus case, introduce the following special complex structure in
${\bf R}^4$, for which $x_1 + ix_2$ and $y_1 - iy_2$ are holomorphic
coordinates. Then $f_2 - if_1 = (x_1 + ix_2) (y_1 -iy_2)$ and the
symplectic form $\omega = \sum dx_i \wedge dy_i$ is the real part of
$d(x_1 + ix_2) \wedge d(y_1 - iy_2)$. After that the proof for the
focus-focus case reduces to the proof for the hyperbolic case, with
complex numbers instead of real numbers.
$\blacksquare$

The local structure of nondegenerate singular points of IHS has been known
for some time. The main result is that locally near a nondegenerate
singular point, the linearized and the nonlinearized systems are the same.
More precisely, if for each IHS we call the singular foliation given by the
connected level sets of the moment map the {\sc associated singular
Lagrangian foliation} (it is a regular Lagrangian foliation outside singular
points of the moment map), then we have:

\begin{thm}[Vey-Eliasson]
\label{thm:local}
Locally near a nondegenerate singular point of an IHS, the associated singular
Lagrangian foliation is diffeomorphic
to that one given by the linearized system. In case of fixed points, i.e.
singular points of maximal corank, the word ``diffeomorphic'' can be
replaced by ``symplectomorphic''.
\end{thm}

The above theorem was proved by
R\"ussmann \cite{Russmann} for the analytic case with $n=2$,
Vey \cite{Vey} for the analytic case with any $n$, Eliasson \cite{Eliasson}
for the smooth case. What R\"ussmann and Vey proved is essentially that
Birkhoff formal canonical transformation converges if the system is
integrable. Note that in Eliasson's
paper there are some minor inaccuracies (e.g. Lemma 6),
and there is only a sketch of the proof
of symplectomorphism for cases other than elliptic. For elliptic case, a new
proof is given by Dufour and Molino \cite{DM}.
Lerman and Umanskii \cite{LU1}
also considered the smooth case with $n=2$, but their results are somewhat
weaker than the above theorem.
Rigorously speaking, the cited papers proved the above theorem only for the
case of fixed points. But the
general case follows easily from the fixed point case.

Thus the local topological structure of the singular Lagrangian foliation
near a nondegenerate singular point is the same as of a linear Poisson action, up
to diffeomorphisms. It is diffeomorphic to a direct product
of a non-singular component and simple singular components. Here a
non-singular component means a local regular foliation of dimension $n-k$
and codimension $n-k$, a simple singular component is a singular
foliation associated to a simple singular fixed point. In particular, if $x$
is a nondegenerate singular point, then we can describe the local level set
of the moment map which contains $x$ as follows: If $x$ is a codimension 1
elliptic fixed point $(n=1)$, then this set is just one point $x$. If $x$ is
a codimension 1 hyperbolic fixed point $(n=1)$, then this set is a ``cross''
(union of the local stable and nonstable curves). If $x$ is a
codimension 2 focus-focus fixed point $(n=2)$ then this set is a
``2-dimensional cross'', i.e. a union of 2 transversally intersecting
2-dimensional Lagrangian planes. In general, the local level set is
diffeomorphic to a direct product of local level sets corresponding to the
above simple cases, and a disk of dimension equal to the rank of the
moment map at $x$. Note that the cross and 2-dimensional cross have a
natural stratification (into one 0-dimensional stratum and four
1-dimensional strata; one 0-dimensional stratum and two
2-dimensional strata respectively). Thus by taking product, it follows
that the local level set of $x$ (i.e. containing $x$) also has a unique
natural stratification. The proof of the following proposition is obvious:

\begin{prop}
\label{prop:locallevelset}
For the natural stratification of the local level set of the moment map at a
nondegenerate singular point $x$ we have: \\
a) The stratification is a direct product of simplest stratifications
discussed above. \\
b) Each local stratum is also a domain in an orbit of the IHS. \\
c) The stratum through $x$ is the local orbit through $x$ and has the
smallest dimension among all local strata. \\
d) If $y$ is another point in this local level set, then $k_e (y) = k_e (x),
k_h (y) \leq k_h(x), k_f (y) \leq k_f (x)$, and the stratum through $y$ has
dimension greater than the dimension of the stratum through $x$ by
$k_h (x) - k_h (y) + 2 k_f(x) - 2 k_f(y)$. \\
e) There are points $y^0 =  y, y^1,..., y^s = x$, $s = k_h (x) - k_h(y) +
k_f (x) - k_f(y)$, in the local level set, such that the stratum of $y^i$
contains the stratum of $y^{i+1}$ in its closure, and
$k_h(y^{i+1}) - k_h (y^{i}) + k_f(y^{i+1}) - k_f (y^{i}) = 1$.
\end{prop}

\section{Associated singular Lagrangian foliation}
\label{section:singleaf}
To understand the topology of  an integrable system, it is
necessary to consider not only local structure of singularities, as in the
previous section, but also
non-local structure of them. By this I mean the structure near the
singular leaves of the Lagrangian foliation associated to an IHS.
We begin by recalling the following definition.

\begin{define}
Let ${\bf F} : (M^{2n},\omega) \to {\bf R}^n$ be the moment map of a given
IHS,  and assume that the preimage  of every point in
${\bf R}^n$ under ${\bf F}$ is compact and the differential $D{\bf F}$ is
nondegenerate almost everywhere.\\
a)  The {\sc leaf  through a point $x_0 \in M$ }
is the minimal closed invariant subset of $M$ which contains $x_0$ and
which does not intersect
the closure of any orbit of the action, except the orbits contained in it.\\
b) From a) it is clear that every point is contained in exactly one leaf.
The (singular) foliation given by these leaves is called the {\sc Lagrangian
foliation associated to the IHS}. The orbit space of the
associated Lagrangian foliation (i.e. the topological space whose points
correspond to the leaves of the foliation, and open sets correspond to
saturated open sets) is called the {\sc orbit space} of the IHS.\\
c) A {\sc singular leaf} is a leaf that contains a singular
point.  An {\sc nondegenerate singular leaf}  is a singular leaf whose
singular points are all nondegenerate.  A singular leaf is called {\sc of
codimension $k$} if $k$ is the maximal corank of its singular points.
\end{define}

According to Arnold-Liouville theorem,  non-singular orbits are Lagrangian
tori of dimension $n$ - hence the notion of associated Lagrangian foliation.
For singular leaves, we use the term codimension
instead of corank since they are nonlinear objects.

\begin{prop}
\label{prop:various}
a) The moment map
${\bf F}$ is constant on every leaf of the associated Lagrangian foliation.\\
b) If the singular point $x_0$ has corank $k$, then the
orbit ${\mathcal O}(x_0)$ of
the Poisson action through it has dimension $n-k$. \\
c) If the closurr of ${\mathcal O}(x_0)$ contains a nondegenerate
singular point $y$, and $y$ does not belong to ${\mathcal O}(x_0)$,
then the corank of $y$ is greater than the corank of $x_0$. \\
d) If the point $x_0$ is nondegenerate, then there is only finitely many
orbits whose closure contains $x_0$. \\
e) In a nondegenerate leaf the closure of every orbit contains only
itself and  orbits of  smaller dimensions. \\
f) Every nondegenerate singular leaf contains only finitely many orbits. \\
g) The orbits of dimension $n-k$ (smallest dimension) in a nondegenerate
leaf of codimension $k$ are diffeomorphic to the $(n-k)$-dimensional torus
${\bf T}^{n-k}$.
\end{prop}

{\it Proof}. Assertion a) follows from the definitions.
To prove b),  notice that,  since the Poisson action
preserves the symplectic structure and the moment map, every point in the
orbit ${\mathcal O}(x_0)$ has the same corank $k$.
We can assume, for example, $dF_1 \wedge dF_2 \wedge ... \wedge dF_{n-k}
\neq 0$ at $x_0$. Then the vector fields  $X_{F_1},...,X_{F_{n-k}}$ generate
a $(n-k)$-dimensional submanifold through $x_0$. Because of the corank $k$,
the other vector field
$X_{F_{n - k +1}}, ..., X_{F_n}$ are tangent to this submanifold, so this
submanifold coincides with ${\mathcal O}(x_0)$. b) is proved.
Assertions c) and d) follow from Proposition \ref{prop:locallevelset}
e) is a consequence of b). f) follows from d) and the compactness assumption.
It follows from e) that the orbits of smallest dimension in $N$
must be compact. Note that, in view
of b), these orbits are orbits of a locally-free ${\bf R}^{n-k}$ action.
Hence they are tori, and g) is proved.
$\blacksquare$

In particular, a singular leaf of codimension $k$ has dimension greater or
equal to $n-k$: it must contain a compact orbit of dimension $n-k$ but it
may also contain non-compact orbits of dimension greater than $n-k$.
If a leaf is nonsingular or
nondegenerate singular, then it is a connected component of the preimage of a
point under the moment map, and a tubular neighborhood of it can be made
saturated with respect to the foliation (i.e. it consists of the whole leaves
only).  Also, the local structure of nondegenerate singular leaves is given
by Proposition \ref{prop:locallevelset}. One can imagine
a singular leaf, at least in case it is nondegenerate, as a stratified
manifold - the stratification given by the orbits of various dimensions,
and the codimension of the leaf equals $n$ minus the smallest dimension of
the orbits.

Later on, a {\sc tubular neighborhood ${\mathcal U}(N)$
of a nondegenerate leaf $N$} will always mean an appropriately
chosen sufficiently small saturated tubular neighborhood.
The Lagrangian foliation
in a tubular neighborhood of a nondegenerate singular leaf $N$ will be
denoted by $({\mathcal U}(N), {\mathcal L})$. Remark that $N$ is a deformation
retract of ${\mathcal U}(N)$. Indeed, one can retract from
${\mathcal U}(N)$ to $N$ by using a appropriate gradient vector field. Thus,
throughout this paper, in topological arguments,
where only the homotopy type matters, we
can replace $N$ by ${\mathcal U}(N)$ and vice versa.

By a {\sc singularity of an IHS} we will mean
(a germ of) a singular foliation $({\mathcal U}(N),{\mathcal L})$.
Two singularities are called {\sc
topologically equivalent} if there is a foliation preserving homeomorphism
(not necessarily symplectomorphism) between them. In this paper all
singularities are assumed to be nondegenerate.

Recall that, given an IHS in $(M^{2n}, \omega)$, its orbit space has a
natural topology which we will use:
the preimages of open subsets of the orbit space are saturated open subsets
of $M^{2n}$. It is the strongest topology for which the projection map from
$M^{2n}$ to the orbit space is continuous.

\begin{prop}
If all singularities of an IHS are nondegenerate, then the corresponding
orbit space is a Hausdorff space.
\end{prop}

{\it Proof}. If all singularities are nondegenerate, then each leaf of the
Lagrangian foliation is a connected level set of the moment map. The moment
map can be factored through the orbit space (denoted by $O^n$):
$$M^{2n} \to O^n \to {\bf R}^n$$
From this the proposition follows easily. $\blacksquare$

Of course, the converse statement to the above proposition is not true.
For the orbit space to be
Hausdorff, it is not necessary that all singularities are nondegenerate. In
fact, IHS's met in mechanics sometimes admit so-called simply-degenerate
singularities, and the orbit space is still Hausdorff.

Thus, the orbit space, under the condition of nondegeneracy, is a good
topological space. A point in the orbit space $O^n$ is called {\sc singular}
if it corresponds to a singular leaf. By Arnold-Liouville theorem, outside
singular points the orbit space admits a unique natural integral affine
structure. Later, from the results of Section \ref{section:action},
one will see that under the additional assumption of topological stability,
the orbit space is a
stratified affine manifold in a very natural sense.
This assumption is called topological stability. Note that this stability is
distinct from the stability in the sense of differential equations, which for
nondegenerate singularities of IHS's is equivalent to the ellipticity
condition.

Let $N$ be a nondegenerate singular leaf. Then $N$ consists of a finitely
many orbits of our IHS. It is well-known that each orbit is of the
type ${\bf T}^{c} \times {\bf
R}^o$ for some nonnegative integers $c,o$.
We would like to know more about these numbers, in
order to understand the topology of $N$.
Recall that $c + o$ is the dimension of the orbit, which is equal to the
rank of the moment map on it.

\begin{define}
Let $N$ be a nondegenerate singular leaf, $x \in N$. If the orbit
${\mathcal O}(x)$ is diffeomorphic to ${\bf T}^c \times {\bf R}^o$, then the
numbers $c$ and $o$ are called {\sc degree of closedness} and {\sc degree of
openness} of ${\mathcal O}(x)$
respectively. The 5-tuple $(k_e,k_h,k_f,c,o)$ is called the {\sc orbit type}
of ${\mathcal O}(x)$, where $(k_e,k_h,k_f)$ is the Williamson type of $x$.
\end{define}

It is clear that the orbit type of an orbit ${\mathcal O}$ does not depend on the
choice of a point $x$ in it. The orbit type of a regular Lagrangian torus is
$(0,0,0,n,0)$ where $n$ is the number of degrees of freedom. In general, we
have $k_e + k_h + 2 k_f + c + o = n$. Moreover:

\begin{prop}
\label{prop:3values}
The following three values: $k_e, k_f + c, k_f + k_h + o$ are invariants of a
singular leaf $N$ of the IHS.
In other words, they don't depend on the choice of an orbit in $N$.
They will be called {\sc degree of ellipticity, closedness, and hyperbolicity}
of $N$, respectively.
\end{prop}

{\it Proof}.
Let $x$ be an arbitrary point in $N$ and $y$ be a
point sufficiently near $x$, not belonging to
${\mathcal O}(x)$. Since $N$ is closed, it suffices to prove that
$k_e (x) = k_e (y) , k_f (x) + c (x) = k_f (y) + c (y)$ and
$k_f (x) + k_h(x) + o(x) = k_f (x) + k_h(y) + o(y)$.

Since $y$ is near to $x$ and $N$ is nondegenerate, it follows that
the closure of ${\mathcal O}(y)$ contains ${\mathcal O}(x)$ and $y$ belongs to the
local level set of the moment map through $x$. By Proposition
\ref{prop:locallevelset}, we can assume that $k_h (y) = k_h (x) - 1,
k_f (y) = k_f (x)$ or $k_h (y) = k_h (x), k_f (y) = k_f (x) - 1.$

We will prove for the case $k_h (y) = k_h (x) - 1,
k_f (y) = k_f (x)$. The other case is similar. The fact $k_e(y) = k_e (x)$
follows from Proposition \ref{prop:locallevelset}. It remains to prove that
$c(y) = c(x)$, since $c(y) + o(y) = c(x) + o(x) + 1$ by Proposition
\ref{prop:locallevelset}. On one hand,
${\mathcal O}(x) = {\bf T}^{c(x)} \times {\bf R}^{o(x)}$ implies that there is
a free ${\bf T}^c(x)$ on ${\mathcal O}(x)$, which is part of the Poisson action
(our IHS).
By using the Poincar\'e map (in
$c(x)$ directions), one can approximate
the generators of this part of the Poisson action by nearby generators
(in the Abelian group ${\bf R}^n$) so
that they  give rise to a ${\bf T}^{c(x)}$ free action on ${\mathcal O}(y)$.
It follows that $c(y) \geq c(x)$. On the other hand, since ${\mathcal O}(x)$ is
in the closure of ${\mathcal O}(y)$, a set of $c(y)$ generators of the Poisson
action, say $X_1,...,X_{c(y)}$, which generate the free ${\bf T}^{c(y)}$
action on ${\mathcal O}(y)$, also generate a ${\bf T}^{c(y)}$
action on ${\mathcal O}(x)$. If this action on ${\mathcal O}(x)$ is locally free
then $c(x) \geq c(y)$ and we are done. Otherwise we would have that
$\sum a_i X_i (x) = 0$ for some nonzero linear combination
$\sum a_1 X_i$ with real coefficients. But then $\sum a_i X_i$ is
small at $y$ (in some metric), therefor
it cannot give rise to period 1 flow on ${\mathcal O}(y)$, and we come to a
contradiction.
$\blacksquare$

As a consequence of the above proposition, we have:

\begin{prop}
If a point $x$ in a singular leaf $N$ has corank $k$ equal to the
codimension of $N$ (maximal possible)
and has Williamson type $(k_e,k_h,k_f)$, then any other
point $x'$ in $N$ with the same corank  will have the same Williamson type.
We will say also that $N$ has the Williamson type $(k_e,k_h,k_f)$.
\end{prop}

{\it Proof}. $k_e$ is invariant because of the previous proposition. It
suffices to proof that $k_f$ is invariant. But we know that $k_f + c$ is
invariant, and for points of maximal corank we have $c = n - k$. $\blacksquare$

The local structure theorems of Williamson, Vey, and Eliasson suggest that
to study the topology of singular Lagrangian foliations we should study the
following simple cases first: codimension 1 elliptic singularities (i.e.
of Williamson type (1,0,0)), codimension 1 hyperbolic singularities
(Williamson type (0,1,0)), and codimension 2 focus-focus singularities
(Williamson type (0,0,1)). The cases of codimension 1 hyperbolic and
codimension 2 focus-focus singularities are considered in the next sections.
The structure of elliptic singularities
(Williamson type (k,0,0)) has been
known for some time, and we will recall it now. The point is that since
elliptic orbits are stable in the sense of differential equations, each
elliptic singular leaf is just one orbit and the non-local problem is just a
parametrized version of the local problem. In particular, for elliptic
singularities the Lagrangian foliations associated to linearized and
nonlinearized systems are the same. To be more precise:

\begin{define}
\label{define:elliptic}
A nondegenerate singular point (singularity) of corank (codimension)
$k$ of an IHS is called {\sc elliptic} if it has Williamson type $(k, 0, 0)$,
i.e. if it has only elliptic components.
\end{define}

\begin{prop}
\label{prop:ell=torus}
A nondegenerate singular leaf of codimension $k$ is elliptic
if and only if it has dimension $n-k$. In this case it is an isotropic
$(n-k)$-dimensional smooth torus, called
{\sc elliptic torus of codimension $k$}.
\end{prop}

\begin{thm}[\cite{Eliasson,DM}]
\label{thm:eliasson}
Let $N$ be an elliptic singularity of codimension $k$ of an IHS
given by the
moment map ${\bf F}: M \to {\bf R}^n$. Then in a tubular
neighborhood of $N$ there is a symplectic system of coordinates \\
$(x_1,\ldots,x_n, y_1  (mod 1),\ldots,y_{n-k} (mod 1), y_{n-k+1},\ldots,
y_n)$, (i.e. $\omega = \sum dx_i \wedge dy_i$), for which
$N = \{ x_1 = \ldots =
x_n = y_{n-k+1} = \ldots = y_n = 0 \}$, such that ${\bf F}$ can be expressed
as a smooth
function of $x_1,\ldots, x_{n-k}, x_{n-k+1}^2 + y_{n-k+1}^2,\ldots, x_n^2 +
y_n^2$:  $$ {\bf F} = {\bf F}(x_1,\ldots, x_{n-k}, x_{n-k+1}^2 +
y_{n-k+1}^2,\ldots, x_n^2 + y_n^2) $$ In other words, the Lagrangian
foliation given by an IHS near every elliptic singularity is
symplectomorphic to the Lagrangian foliation associated to the linearized
system.
\end{thm}

See the paper of Dufour and Molino \cite{DM} for a proof of the above
theorem. In \cite{Eliasson}, Eliasson simply says that it is just a
parametrized version of his local structure theorem.

\begin{cor}
If $N$ is an elliptic singular leaf of codimension $k$ of an IHS with $n$
degrees of freedom then there is a unique natural $T^n$ symplectic action
in ${\mathcal U}(N)$ which preserves the moment map, and which is free almost
everywhere. There is also a highly
non-unique ${\bf T}^{n-k}$ subgroup action of the above action, which is
free everywhere in ${\mathcal U}(N)$.
\end{cor}

Using one of the free ${\bf T}^{n-k}$ actions in the above corollary, one
can use the Marsden-Weinstein reduction
procedure (cf. \cite{MW}) to reduce codimension $k$ elliptic
singularities of IHS's with $n$ degrees of freedom to codimension $k$
elliptic singularities of IHS's with $k$ degrees of freedom.

\section{Codimension 1 case}
\label{section:codim1}
All known integrable systems from physics and mechanics exhibit codimension 1
nondegenerate singularities. For these singularities, it follows from Williamson's
classification that there are only 2 possible cases: elliptic and hyperbolic.
Elliptic singularities were discussed in the previous section. In this section we
assume all singularities to be hyperbolic, though the results will be also obviously
true for elliptic singularities (except the uniqueness of the torus ${\bf T}^{n-1}$
action).

\begin{thm}
\label{thm:Tn-1action}
Suppose $N$ is a nondegenerate codimension 1 singularity of an IHS with $n$
degrees of freedom. Then in $\mathcal{U}(N)$ there is a Hamiltonian
${\bf T}^{n-1}$ action such that: \\
a) This group action preserves the moment map. \\
b) This group action is locally free and it
is free outside singular points of $\mathcal{U}(N)$. \\
c) Each singular point can have at most one non-trivial element of
${\bf T}^{n-1}$ which preserves it. In other words, the isotropy group of
each point is at most ${\bf Z}_2$.
\end{thm}

{\bf Remark}. In case $n=2$, the above theorem appeared in \cite{Zung90}.

{\it Proof}.
Since $N$ is hyperbolic, it contains not only singular points, but also
regular points. Let $y$ be a regular point in $N$, then the orbit
${\mathcal O}(y)$ will be diffeomorphic to ${\bf T}^{n-1} \times {\bf R}^1$,
according to the previous section. We can choose $n-1$ generators of the
Poisson action of our IHS, in terms of $n-1$ Hamiltonian vector fields
$X_1,...,X_{n-1}$, such that the flow of each $X_i$ restricted to ${\mathcal O}(y)$
is periodic with period 1, and together they generate a free ${\bf T}^{n-1}$
action on ${\mathcal O}(y)$. Let ${\mathcal  O}(x)$ be a singular orbit contained in
the closure of ${\mathcal O}(y)$. Then of course the flow of $X_i$ are also time
1 periodic in ${\mathcal O}(x)$. Let ${\mathcal O}(z)$ be another regular orbit in
$N$ which also contains ${\mathcal O}(x)$. Let $y' \in {\mathcal O}(y)$,
$z' \in {\mathcal O}(z)$ are two points sufficiently close to $x$. We will
assume that they lie in different n-disks of the cross (times a disk)
in the local level set at $x$ (see Proposition
\ref{prop:locallevelset}). Then it
follows from the local structure of the Lagrangian foliation at $x$ that
there is a regular orbit ${\mathcal O}_1$ (not contained in $N$) passing
arbitrarily near to both points $y$ and $z$. It follows that there are
generators $X_i'$ near to $X_i$ of our  Poisson action
such that the flows of $X_i$ are periodic of
period exactly 1 in ${\mathcal O}_1$. In turn, there are generators $X_i''$
near to $X_i'$ such that the flows of $X_i''$ are periodic of period exactly
1 in ${\mathcal O}(z)$. By making ${\mathcal O}_1$ tend to pass through $y$ and $z$,
in the limit we obtain that the flows of $X_i$ are periodic of period
exactly 1 in ${\mathcal O}(z)$. Since $N$ is connected, by induction we see that
$X_1,...X_{n-1}$ generates a ${\bf T}^{n-1}$ action on $N$, and this action
is free on regular orbits (of type ${\bf T}^{n-1} \times {\bf R}$) in $N$.
From the proof of Proposition \ref{prop:3values} it follows that this action
is also locally free on singular orbits of $N$.

To extend the above ${\bf T}^{n-1}$ action from $N$ to ${\mathcal U}(N)$, one can
use the same process of extending the vector fields from ${\mathcal O}(y)$ to
${\mathcal O}_1$ as before. As a result, we obtain a natural
${\bf T}^{n-1}$ action in ${\mathcal U}(N)$, unique up to automorphisms of
${\bf T}^{n-1}$.
Obviously, this action preserves our IHS.
Outside singular leaves, this action is Hamiltonian,
because of Arnold-Liouville theorem. Since singular leaves in ${\mathcal U}(N)$
form a small subset of measure 0, it follows that the above action is
symplectic in the whole ${\mathcal U}(N)$. Since the symplectic form
$\omega$ is exact in ${\mathcal U}(N)$ (because it is obviously exact in $N$),
this torus action is even Hamiltonian.

We have proved assertions a) and b). To prove c), assume that $\xi$ is a
nonzero linear combination of $X_i$, which generates a periodic flow in
${\mathcal O}(x)$ of period 1, and denote by $g_1$ its time 1 map in $N$. Since
$y'$ is near to $x$, $g_1 (y')$ is also near to $x$. If $g_1(y')$ belongs to
the same local stratum at $x$ as $y'$, then it follows that $g_1(y')$ must
coincide with $y'$, and $\xi$ is an integral linear combination of $X_i$. If
not, $g_1(y')$ must line in the local stratum ``opposite'' to the
local stratum containing $y'$, and it will follow that twice of $\xi$ is an
integral combination of $X_i$. Thus, any isotropic element of the ${\bf
T}^{n-1}$ action on ${\mathcal O}(x)$ has order at most 2. To see that the
isotropic subgroup has at most two elements, assume for example
that  $X_1/2$ and
$X_2/2$ generate periodic flows of period 1 in ${\mathcal O}(x)$. Then the time
1 maps of these vector fields send $y'$ to the opposite local stratum. It
follows that the time 1 map of $(X_1 + X_2)/2$ sends $y'$ to itself, and
$(X_1 + X_2)/2$ is an integral linear combination of $X_i$, a contradiction.
$\blacksquare$

\begin{prop}
\label{prop:doublecover} If in the above theorem, there are points with non-trivial
isotropy group, then there is a natural unique\footnote{erratum: the double covering
is not always unique, though there are only a finite number of possibilities} double
covering of ${\mathcal U}(N)$, denoted by $\overline{{\mathcal U}(N)}$, such that
everything (symplectic form, IHS, ${\bf T}^{n-1}$-action) can be lifted to this
double covering, and the ${\bf T}^{n-1}$-action in this covering will be free.
\end{prop}

{\it Proof}. The above theorem, without the torus action, was in fact proved
by Fomenko \cite{Fomenko}. To prove it, one can write down a presentation
for the fundamental group $\pi_1(N) = \pi_1({\mathcal U}(N))$ as follows:

{\it Generators} (three types): $\alpha_1,...,\alpha_{n-1}$, which are generated by
$X_1,...,X_{n-1}$ in a regular orbit in $N$. $\beta_1,...,\beta_s (s \geq 0)$, which
are ``exceptional cycles'': they lie in singular orbits and are not conjugate to an
integral combination of $\alpha_i$. $\gamma_1,...,\gamma_t (t \geq 1)$, which are
``base cycles''\footnote{erratum: as regards $\gamma_j$, they are not unique even
when their projections to the orbit space of the ${\bf T}^{n-1}$-action are fixed:
for example, $\gamma_1$ may be changed to $\gamma_1\beta_1$, and that explains why
the double covering exists but is not unique in general}: the fundamental group of
the quotient space of $N$ by the ${\bf T}^{n-1}$ action is generated by the image of
these cycles.

{\it Relations}: $\alpha_i$ commute with all the other generators. Twice of
$\beta_i$ are integral combinations of $\alpha_i$.

Let $G$ be the subgroup of $\pi_1(N)$, consisting of the words which contain
an even total number of $\beta_i$. For example, $\alpha_1$ and
$\beta_1 \beta_2$ are elements of $G$. Then $G$ is a subgroup of index 2,
and it is easy to see that the double covering associated to $G$ will
satisfy our requirements.
$\blacksquare$

In the following theorem, which gives a canonical form for codimension 1
nondegenerate hyperbolic singularities, we assume that the torus action
discussed above is free. Otherwise we will take the double covering first.

\begin{thm}
\label{thm:canonical1h}
There exist a system of coordinates
$(x_1, y_1, \ldots , x_n, y_n)$ in (an appropriately chosen) ${\mathcal U}(N)$,
where $y_1,\ldots,y_{n-1}$ are defined modulo 1 (angle coordinates), and
$(x_n, y_n)$ defines an immersion from a surface
$P^2$ to ${\bf R}^2$, such that: \\
a) These coordinates give a diffeomorphism from  ${\mathcal U}(N)$
to ${\bf D}^{n-1} \times {\bf T}^{n-1} \times P^2$ \\
b) The symplectic form $\omega$ has the canonical form
$$\omega = \sum_{i = 1}^{n-1} dx_i \wedge dy_i + \pi^{\ast}(\omega_1)$$
where $\omega_1$ is some area form on $P^2$, and $\pi^{\ast}$ means the
lifting. \\
c) $x_1,\ldots, x_{n-1}$ are invariants of our IHS. \\
Such a system of coordinates is not unique, but will be called
{\sc canonical}.
\end{thm}

{\it Proof}. Let $x_i$ be the Hamiltonian functions of Hamiltonian vector
fields which generate the ${\bf T}^{n-1}$ action
in Theorem~\ref{thm:Tn-1action}, $x_i(N) = 0$. Denote their corresponding
vector fields by $\xi_i$.
Remark that the ${\bf T}^{n-1}$-group action in Theorem~\ref{thm:Tn-1action}
gives rise to a trivial ${\bf T}^{n-1}$-foliation.
Denote its base by $B^{n+1}$. Let $L$ be a section of this foliation, and
define functions $z_i \; (i = 1,\ldots,n-1)$
by putting them equal to zero on $L$ and setting
$dz_i(\xi_i) = \{x_i, z_i\} = 1$. Set
$\omega_1 = \omega - \sum_1^{n-1} dx_i \wedge dz_i$.
Then one checks that $L_{\xi_i} \omega_1 = i_{\xi_i} \omega_1 = 0$. It means
that $\omega_1$ is a lift of some closed 2-form from $B^{n+1}$ to
${\mathcal U}(N)$, which we will also denote by $\omega_1$. Since $\omega$
is non-degenerate, it follows that $\omega_1$ is non-degenerate on every
2-surface (with boundary) $P^2_{x_1,\ldots,x_{n-1}} =
B^{n+1} \cap \{x_1,\ldots,x_{n-1} {\tt fixed}\}$. Using Moser's path method
\cite{Moser}, one can construct a diffeomorphism
$\phi : B^{n+1} \to {\bf D}^{n-1} \times P^2$, under which $\omega_1$
restricted to $P^2_{x_1,\ldots,x_{n-1}}$ does not depend on the choice of
$x_1,\ldots,x_{n-1}$. In other words, there is an area form $\omega_2$
on $P^2$ such that $\omega_1 - \omega_2$ vanishes on every
$P^2_{x_1,\ldots,x_{n-1}}$. Since $d(\omega_1 - \omega_2) = 0$, we can write
it as $\omega_1 - \omega_2 = d(\sum_1^{n-1}
a_i dx_i + \beta)$, where $\beta$ is
some 1-form on $B^{n+1}$ (which is not zero on $P^2_{x_1,\ldots,x_{n-1}}$
in general). If we can eliminate
$\beta$, i.e. write $\omega_1 - \omega_2 = d(\sum_1^{n-1} a_i dx_i)$,
then we will have $\omega = (\sum_1^{n-1}
dx_i \wedge dz_i - \sum_1^{n-1} dx_i \wedge da_i + \omega_2
= \sum_1^{n-1}dx_i \wedge d(z_i - a_i) + \omega_2$, and the theorem will
be proved by  putting $y_i = z_i - a_i$. Let us show now how to
eliminate $\beta$. $\beta$ restricted on every $P^2_{x_1,\ldots,x_{n-1}}$
is a closed 1-form, hence it represents a cohomology element
$[\beta](x_1,\ldots,x_{n-1}) \in H^1(P^2)$. If
$[\beta](x_1,\ldots,x_{n-1}) \equiv 0$ then
$\beta = dF - b_1 dx_1 - \ldots - b_{n-1} dx_{n-1}$
for some functions $F, b_1,\ldots,b_n$, and we have
$\omega_1 - \omega_2 = d(\sum(a_i - b_i x_i) dx_i)$. In general, we can
achieve $[\beta](x_1,\ldots,x_{n-1}) \equiv 0$ by induction on the number of
generators of $H^1(P^2)$ as follows. Let $\gamma$ be a simple curve in $P^2$
which represents a non-zero cycle. Set
$b(x_1,\ldots,x_{n-1}) = <[\beta], \gamma>(x_1,\ldots,x_{n-1})$.
Immerse $P^2$ in an annulus so that only simple curves homotopic to $\gamma$
go to non-zero cycle there, other simple curves go to vanishing cycles.
By this immersion we have a (non single-valued) system of coordinates
$(u,v)$ on $P^2$, $u - {\rm mod}\; 1$. $\omega_2$ has the form
$\omega_2 = a du \wedge dv$ for some positive function $a$. Change
$\omega_2$ for the following 2-form on ${\mathcal U}(N)$:
$\omega_2' = \omega_2 + db \wedge du$. It is clear that $\omega_2'$ and
$\omega_2$ restricted on every $P^2_{x_1,\ldots,x_{n-1}}$ are the same.
Moreover, $\omega_2'$ is closed and of rank 2. Thus the distribution by its
$(n-1)$-dimensional tangent zero-subspaces is integrable, and it gives rise
again to a diffeomorphism $\phi' : B^{n+1} \to {\bf D}^{n-1} \times P^2$.
Replacing $\omega_2$ by $\omega_2'$, we have
$\omega_1 - \omega_2' = d(\sum_1^{n-1}
a_i dx_i + \beta')$, with $\beta' = \beta - b du$, whence
$<[\beta'],\gamma> = 0$. $\blacksquare$

Using the ${\bf T}^{n-1}$ action constructed in the previous theorems, we
can apply the Marsden-Weinstein reduction procedure to obtain a
$(n-1)$-dimensional family of codimension 1 hyperbolic singularities of IHS
with one degree of freedom. These singularities are called {\sc surface
singularities} because they lie in 2-dimensional surfaces.
They are so simple that we are
not going to make them simpler. In \cite{BFM}, Fomenko called them
``letter-atoms''. These ``letter-atoms'' were computed for many interesting
examples of IHS arising in classical mechanics and physics (see e.g.
\cite{AS,Fomenko2,GOC,Kharlamov,Oshemkov,ZungLada}).
Some of these codimension 1 singularities even have special names (see
\cite{BFM}).

After the Marsden-Weinstein reduction we have an $(n-1)$-dimensional family
of surface singularities. There is no guarantee that surface singularities in
this family must be topologically equivalent, although in all known
examples of IHS arising in mechanics
and physics they turn out to be so. This situation leads to the following:

\begin{prop}
\label{prop:stable1h}
For a codimension 1 hyperbolic singularity ${\mathcal U}(N)$
of an IHS with $n$ degrees of
freedom, the following three conditions are equivalent: \\
a) All singular leaves in ${\mathcal U}(N)$ are topologically equivalent. \\
b) Under the Marsden-Weinstein reduction with respect to the
torus action given by Theorem \ref{thm:Tn-1action},
the topological structure of
surface singularities in the obtained \mbox{$(n-1)$}-dimensional
family is constant (i.e. does not depend on the parameter). \\
c) The singular value set of the moment map restricted to ${\mathcal U}(N)$ is
a smooth $(n-1)$-dimensional disk containing the image of $N$ in
${\bf R}^n$.
\end{prop}

\begin{define}
\label{define:stable-codim1}
A non-degenerate codimension 1 singularity of an IHS is
called {\sc topologically stable} if it satisfies one of the equivalent
conditions in Proposition \ref{prop:stable1h}.
\end{define}

{\it Proof} of Proposition \ref{prop:stable1h}. The proof is rather
straightforward. We will prove that c) implies a).
Without loss of generality, we can assume that the local
singular value set in ${\bf R}^n$ is given by the equation $F_1 = 0$, where
$F_1$ is one of the components of the moment map. Then $F_1 = 0$ on all
singular leaves in ${\mathcal U}(N)$. Let ${\mathcal O}_1$ and ${\mathcal O}_2$ be two
hyperbolic orbits in $N$. It follows from the nondegeneracy that $F_1 (y)
\neq 0$ for any point $y$ near to ${\mathcal O}_1$ and not belonging to a local
singular leaf near to ${\mathcal O}_1$. Since by continuity,
any singular leaf near to ${\mathcal O}_2$ passes nearby ${\mathcal O}_1$, it follows
that each such a singular leaf must contain a local singular leaf
near to ${\mathcal O}_1$. This fact implies the topological equivalence of all
singular leaves near to $N$.
$\blacksquare$

{\bf Remark}. In case $n=2$ the notion of topological stability (for
codimension 1 hyperbolic singularities) appeared in \cite{Bolsinov}.

Theorem \ref{thm:canonical1h} tells us that if the ${\bf T}^{n-1}$ action on
$\mathcal{U}(N)$ is free then $\mathcal{U}(N)$ is symplectomorphic to
${\bf D}^{n-1} \times {\bf T}^{n-1} \times P^2$. If the action is not free
then $\mathcal{U}(N)$ is symplectomorphic to the quotient of the above product
by a free action of some finite group $\Gamma$. This group $\Gamma$ may be
taken to be ${\bf Z}_2$, according to Proposition \ref{prop:doublecover}.
However, if we require that $\Gamma$ act on
${\bf D}^{n-1} \times {\bf T}^{n-1} \times P^2$
component-wise (i.e. the action commutes with the projections onto
the components of the product, then in general it must be bigger than
${\bf Z}_2$. Indeed, if ${\bf Z}_2$ acts on
$\overline{{\mathcal U}(N)} = {\bf D}^{n-1} \times {\bf T}^{n-1} \times P^2$
with quotient equal to ${\mathcal U}(N)$, then exceptional cycles
$\beta_1,...\beta_s$ of $\pi_1(N)$
(as in Proposition \ref{prop:doublecover})
must have the same image in
${\bf T}^{n-1}$. It is clearly not the case in
general (except if $n=2$), so we have to modify our construction of finite
covering of ${\mathcal U}(N)$ as follows:

Remember that we have generators $\alpha_1,..., \alpha_{n-1}; \beta_1,..., \beta_s;
\gamma_1,..., \gamma_t$ for $\pi_1(N)$ (cf. Proposition \ref{prop:doublecover}). An
important observation is that there is a natural homomorphism from $\pi_1(N)$ to the
group ${\bf T}^{n-1}$, associated  to the action of ${\bf T}^{n-1}$ on $N$, which
sends $\alpha_1,..., \alpha_{n-1}, \gamma_1,..., \gamma_t$ to 0, and $\beta_1,...,
\beta_s$ to elements of order 2. Let $G_{can} \supset \pi_1(N)$ be the kernel of
this homomorphism, and $\Gamma_{can} = \pi_1(N) / G_{can}$ be the image. So
$\Gamma_{can}$ is a subgroup of ${\bf T}^{n-1}$ (which contains only elements of
order at most 2). Denote by $\overline{{\mathcal U}(N)}_{can}$ the normal finite
covering of ${\mathcal U}(N)$ corresponding to $G_{can}$, and call it the canonical
covering. Then we have:

\begin{thm}
\label{thm:AL-1h}
Let $({\mathcal U}(N), {\mathcal L})$ be a nondegenerate
codimension 1 singularity of an IHS with $n$ degrees of freedom. \\
i) The symplectic form, IHS, and ${\bf T}^{n-1}$-action can be lifted from
${\mathcal U}(N)$ to the canonical covering $\overline{{\mathcal U}(N)}_{can}$,
and the ${\bf T}^{n-1}$-action on this canonical covering will be free. \\
ii)  $\overline{{\mathcal U}(N)}_{can}$ is symplectomorphic to
${\bf D}^{n-1} \times {\bf T}^{n-1} \times P^2$
with the canonical symplectic form as in Theorem \ref{thm:canonical1h}.
Under this symplectomorphism, $\Gamma_{can}$ acts on
${\bf D}^{n-1} \times {\bf T}^{n-1} \times P^2$
component-wise, i.e. it commutes with the projections. We say
that $\overline{{\mathcal U}(N)}_{can}$ admits an {\sc equivariant canonical
system of coordinates}.
\end{thm}

{\it Proof}.
i) The proof is immediate from the definition of
$\overline{{\mathcal U}(N)}_{can}$. \\
ii) $\overline{{\mathcal U}(N)}_{can}$ is a trivial principal ${\bf
T}^{n-1}$-bundle, and $\Gamma_{can}$ acts on its sections. By lifting a
section from $N$, we can choose a section $L_1$ in
$\overline{{\mathcal U}(N)}_{can}$
which is stable under $\Gamma_{can}$-action homotopically. ($N$ is a
singular ${\bf T}^{n-1}$ bundle with base space equal to a graph, so we can
talk about its sections). Define
$$ L_2 (x) = 1/ (\# \Gamma_{can}) \sum_{\gamma \in \Gamma_{can}}
( (\gamma L_1)(x) - \gamma )$$
where $x$ is a point in the base space of $\overline{{\mathcal U}(N)}_{can}$,
$\gamma L_1$ is the action of $\gamma$ on $L_1$, $ - \gamma$ means the action
of $ - \gamma$ in the torus over $x$. Then $L_2$ is a section which is
$\Gamma_{can}$-equivariant, that is $\gamma L_2 = L_2 + \gamma$ for any
$\gamma \in \Gamma_{can}$. Now applying arguments of the proof of Theorem
\ref{thm:canonical1h} to $L_2$, keeping the
equivariance of $L_2$ all the time,
we will find an equivariant coisotropic section $L_3$. It is immediate that
if we trivialize $\overline{{\mathcal U}(N)}_{can}$ via $L_3$, $\Gamma_{can}$
will act on this trivialization component-wise.
$\blacksquare$

{\bf Remark}. It is also clear that $\Gamma_{can}$ is the smallest group
that acts component-wise on something of the type
${\bf D}^{n-1} \times {\bf T}^{n-1} \times P^2$ with
${\mathcal U}(N)$ as the quotient. We will call
$(\overline{{\mathcal U}(N)}_{can}, {\rm action \; of} \; \Gamma_{can})$ the
{\sc canonical model} of singularity $({\mathcal U}(N), {\mathcal L})$ because of
this property.

\section{Codimension 2 focus-focus case}
\label{section:focus}
There are quite many integrable systems which exhibit focus-focus
singularities. For example, consider a Garnier system with 2 degrees of
freedom, given by the Hamiltonian $H = 1/2 (p_x^2 + p_y^2) - \alpha_1 x^2
- \alpha_2 y^2 + (x^2 + y^2)^2$, where $0 < \alpha_1 < \alpha_2$. The graph
of the potential $V =  - \alpha_1 x^2
- \alpha_2 y^2 + (x^2 + y^2)^2$ of this system looks like the bottom of a
wine bottle. In particular, the point $x=y = 0$ is a nondegenerate local
maximum of $V$, and it follows immediately that the point $p_x = p_y = x = y
= 0$ is focus-focus in the phase space. Topological classification of
focus-focus singularities for integrable systems with two degrees of freedom
was obtained in \cite{ZungFocus}. Let us recall the main result from there.

\begin{thm}[\cite{ZungFocus}]
\label{thm:focusn=2}
Let $N$ be nondegenerate singular leaf in an IHS with two degrees of freedom
which contains a focus-focus (fixed) point. Then we have: \\
a) There is a unique natural ${\bf S}^1$ Hamiltonian action in a
neighborhood of $N$ which preserves the IHS. This action is locally free
outside focus-focus fixed points of $N$. \\
b) If $N$ does not contains closed 1-dimensional orbits, then the local
orbit space is a disk with a ``removable'' singular point at the image of
$N$ under the projection. However, it is not affinely equivalent to a
regular affine disk. In fact, the monodromy obtained by parallel
transportation around the image of $N$ of the affine structure is
represented by the matrix
$
\begin{pmatrix}
  1 & m \\
  0 & 1
\end{pmatrix},
$
where $m$ is the
number of fixed focus-focus points in $N$. Moreover, any two
such singularities with the same number of focus-focus points are
topologically equivalent. \\
c) If $N$ contains closed 1-dimensional orbits, then there is an arbitrarily
$C^{\infty}$ small integrable perturbation of the IHS, under which $N$ is
replaced by a new singular leaf $N'$ which is close to $N$, contains the
same number of focus-focus points as $N$, and contains no 1-dimensional
closed orbit.
\end{thm}

The above theorem describes the structure of codimension 2 focus-focus singularities
in case of 2 degrees of freedom. Assertion a) can be seen from the local structure
of focus-focus singular points. Assertion b) is called the phenomenon of
non-triviality of the monodromy, which was first observed by Duistermaat and Cushman
\cite{Duistermaat} for the spherical pendulum. This nontriviality of the monodromy
was then found for various systems (see e.g. \cite{CK,Bates}) before we observed
that it is a common property of focus-focus singularities. The proof of b) given in
\cite{ZungFocus} is rather simple and purely topological. Independently, Zou
\cite{Zou} also proved b), for the case $m=1$, under some additional assumptions.
Lerman and Umanskii \cite{LU1,LU3} also studied the topology of focus-focus
singularities, but their description seems too complicated. Assertion c) says that
the condition in b) can be always achieved by a small perturbation. c) can be proved
easily by the use of the ${\bf S}^1$ action given in a).

The condition given in b) of the above theorem is called the {\sc
topological stability} condition. In other words, a focus-focus singular
leaf (i.e. a singular leaf containing focus-focus points) in an IHS with two
degrees of freedom is called topologically stable if it is ``purely
focus-focus'', i.e. if it does not contain codimension 1 hyperbolic singular
points.

Like in previous section, we would like to reduce codimension 2
focus-focus singularities (i.e. singularities which contain codimension 2
focus-focus points) of IHS's with many degrees of freedom to focus-focus
singularities of IHS's with 2 degrees of freedom. First of all, we need some
torus action. This action is provided by the following theorem.

\begin{thm}
\label{thm:focus}
Let $N$ be a codimension 2 nondegenerate
focus-focus singular leaf in an IHS with $n$
degrees of freedom, $n > 2$. Then we have: \\
a) There is a natural ${\bf T}^{n-1}$ action in a tubular saturated
neighborhood ${\mathcal U}(N)$ of $N$, which preserves the given IHS.
This group action is unique up to isomorphisms of the Abelian group
${\bf T}^{n-1}$. \\
b) If ${\mathcal O}$ is an orbit in $N$ consisting of codimension 2 focus-focus
points, then ${\mathcal O}$ is diffeomorphic to ${\bf T}^{n-2}$, and the above
action is transitive on ${\mathcal O}$. Moreover, the isotropy group is ${\bf
S}^1$, i.e. a connected 1-dimensional subgroup of ${\bf T}^{n-1}$. \\
c) If ${\mathcal O}'$ is an orbit in $N$ consisting of codimension 1 hyperbolic
singular points, then the above action is transitive and locally free on
${\mathcal O}'$. Moreover, the isotropy group is at most ${\bf Z}_2$. \\
d) All hyperbolic orbits in $N$ can be taken away by a $C^{\infty}$
small integrable perturbation of the IHS, like in assertion c) of Theorem
\ref{thm:focusn=2}
\end{thm}

{\it Proof}. If $\mathcal O$ is a regular orbit in $N$, then it follows from
Proposition \ref{prop:3values} that $\mathcal O$ is diffeomorphic to ${\bf T}^{n-1}
\times {\bf R}^1$. Consequently, there is a natural ${\bf T}^{n-1}$ action in
${\mathcal O}$. By the same arguments as in Theorem \ref{thm:Tn-1action}, one can
extend this action to a natural Hamiltonian ${\bf T}^{n-1}$ action in ${\mathcal
U}(N)$, which preserves the moment map, and a) is proved. Assertion c) is also
proved by the same arguments as in Theorem \ref{thm:Tn-1action}. From the existence
of torus action, one can apply the canonical coordinates of Theorem
\ref{thm:canonical1h} to easily prove assertion d). Only assertion b) needs a little
more work.

To prove b), let ${\mathcal O}(x)$ be a focus-focus singular orbit in $N$. Then it
is diffeomorphic to ${\bf T}^{n-2}$. Since the ${\bf T}^{n-1}$ action is transitive
in ${\mathcal O}(x)$, the isotropy group of (or at) ${\mathcal O}(x)$ is a closed
1-dimensional subgroup of the Abelian group ${\bf T}^{n-1}$. We have to show that
this subgroup is connected.

We can choose $n-1$ generators of the ${\bf T}^{n-1}$  action in terms of
vector fields $X_1,...X_{n-1}$, such that the flows of these vector fields
are periodic with minimal period 1 in $N$, and the vector field $X_1$
generates the connected component of the isotropy group of ${\mathcal
O}(x)$.Assume by contradiction that there is a non-integral
linear combination $\xi= \sum a_i X_i$, whose flow is time 1 periodic in
${\mathcal O}(x)$. Say $a_2$ is not an integral number. Denote the time 1 map of
$\xi$ by $g$. Of course, $g(x) = x$.
Let $y$ be a point in $N$ which is near to $x$ and not belonging to
${\mathcal O}(x)$. Then it follows from assertion c) of Proposition
\ref{prop:localautogroup} that $g(y)$ lies in the same stratum of the local
level set (of the moment map at $x$) as $y$. Consequently, there is another
linear combination ${\xi}' = \sum a_i' X_i$, with $a_i'$ near to $a_i$ if $i
\neq 1$, such that the time one map of ${\xi}'$ preserves $y$. But then
on one hand, $a_i'$ must be integral, and on the other hand, $a_2'$ is near
to $a_2$ and cannot be integral. This contradiction shows that the isotropy
group at ${\mathcal O}(x)$ of the ${\bf T}^{n-1}$ action must be connected.
$\blacksquare$

We can apply the Marsden-Weinstein reduction procedure to the above
${\bf T}^{n-1}$ action. But notice that, since the action is far from being
free, in the result we will obtain singular symplectic 2-dimensional
orbifolds, but not manifolds. It is what was done in \cite{ZungFocus}
for the case $n=2$. To reduce to focus-focus singularities of systems with 2
degrees of freedom, however, we need an appropriately chosen ${\bf T}^{n-2}$
subgroup action of our ${\bf T}^{n-1}$ group action. Denote such a subgroup in
${\bf T}^{n-1}$ by ${\bf T}^{n-2}_0$. First of all, the action of ${\bf T}^{n-2}_0$ on
singular focus-focus orbits must be at least locally free. It means that
${\bf T}^{n-2}$ is transversal to $S_1,...,S_m$ in ${\bf T}^{n-1}$, where $S_i$ are
isotropy groups at singular focus-focus orbits, denoted by  ${\mathcal O}_i$,
$m$ being the number of such orbits in $N$. Second, to get a regular
reduction, the action of ${\bf T}^{n-2}_0$ on ${\mathcal O}_i$ must be free but not
only locally free. If $N$ does not contain hyperbolic orbits, then
this condition is equivalent to the condition that
${\bf T}^{n-2}_0$ intersects each of $S_i$ only at 0.

Notice that $S_i$ represent vanishing cycles in $N$, and they can kill a big part of
$\pi_1({\bf T}^{n-1})$ in the fundamental group $\pi_1({\mathcal U}(N))$ of
${\mathcal U}(N)$. As a consequence, if such a subgroup ${\bf T}^{n-2}_0$ as above
does not exist, then in general there is no finite covering of ${\mathcal U}(N)$
admitting such a subgroup. In other words, there is no finite covering argument like
in Proposition \ref{prop:doublecover}. However, the following condition of
topological stability will ensure that a subgroup ${\bf T}^{n-2}_0$ of ${\bf
T}^{n-1}$ which intersects $S_i$ only at 0 does exist. In fact, this condition will
imply that all 1-dimensional subgroups $S_i$ coincide.

\begin{define}
\label{define:focusstable}
A codimension 2 nondegenerate focus-focus singularity ${\mathcal U}(N)$
is  called {\sc topologically stable} if it has the following two properties:
\\
a) $N$ does not contain hyperbolic singular orbits. \\
b) All singular leaves in ${\mathcal U}(N)$ are topologically equivalent.
\end{define}

\begin{prop}
\label{prop:stable2f}
Let ${\mathcal U}(N)$ be a codimension 2 nondegenerate focus-focus
singularities. \\
a) If ${\mathcal U}(N)$ satisfies the condition b) in Definition
\ref{define:focusstable} then all of the isotopy groups $S_i$
at focus-focus orbits in $N$ of the symplectic ${\bf T}^{n-1}$ action
(given by Theorem \ref{thm:focus} are the same. \\
b) ${\mathcal U}(N)$ is topologically stable if and only if
the singular value set of the moment map restricted to ${\mathcal U}(N)$ is
a smooth $(n-2)$-dimensional disk containing the image of $N$ in
${\bf R}^n$.
\end{prop}

{\it Proof}. a) Suppose ${\mathcal U}(N)$ satisfies the condition b) of
Definition \ref{define:focusstable}. Let ${\mathcal O}_1$ be a focus-focus orbit
in $N$. Let $X_1,...,X_{n-1}$ be a system of generators of the
${\bf T}^{n-1}$ action as in Theorem \ref{thm:focus}, given in terms of
Hamiltnonian vector fields, such that $X_1$ vanishes in ${\mathcal O}_1$. Then
$X_1$ generates the isotropy group at ${\mathcal O}_1$. Let $f_1$ be the
Hamiltonian function associated to $X_1$, defined by $f_1(N) = 0$. It is
easy to be seen that $X_1$ also generates isotropy groups of focus-focus
orbits near to ${\mathcal O}_1$ (these orbits lie in ${\mathcal U}(N)$ but outside
$N$). Consequently, $X_1$ vanishes on all focus-focus orbits near to
${\mathcal O}_1$, and $f_1 = 0$ on all these orbits. Let ${\mathcal O}_2$ be another
focus-focus orbit in $N$. The condition b) of
Definition \ref{define:focusstable} implies that singular leaves which
contain focus-focus orbits near to ${\mathcal O}_2$ also contain focus-focus
orbits near to ${\mathcal O}_1$. Since $f_1$ is invariant on each leaf, it
follows that $f_1 = 0$ on all focus-focus orbits near to ${\mathcal O}_2$.
Consequently, the symplectic gradient $X_1$ of $f_1$ vanishes
on ${\mathcal O}_2$, and hence $X_1$ generates the isotropic group at ${\mathcal
O}_2$. Thus the isotropic groups at ${\mathcal O}_1$ and ${\mathcal O}_2$ are the
same.

b)  If ${\mathcal U}(N)$ is topologically stable then clearly
the singular value set of the moment map restricted to ${\mathcal U}(N)$ is
a smooth $(n-2)$-dimensional disk containing the image of $N$ in
${\bf R}^n$.  We will prove the inverse statement. Suppose that the local
singular value set of the moment map in ${\bf R}^n$ is a disk of dimension
$(n-2)$. Then $N$ cannot contain hyperbolic codimension 1 singular orbits,
otherwise this local singular value
set must contain a subset diffeomorphic to ${\bf R}^{n-1}$.
The rest of the proof resembles that of Proposition \ref{prop:stable1h}.
$\blacksquare$

Thus, for stable codimension 2 focus-focus singularities, one can use
Marsden-Weinstein reduction with respect to some choice of ${\bf T}^{n-2}_0$
to reduce them to stable focus-focus singularities of IHS's with two degrees
of freedom. Notice that, like the case of elliptic singularities and unlike
the case of codimension 1 hyperbolic singularities, the free ${\bf T}^{n-2}$
exists but is not unique.

\section{Torus action and topological stability}
\label{section:action}
In this section we generalize the results obtained in the previous sections
about torus actions and topological stability to nondegenerate singularities
of any Williamson type.

Recall that if  $({\mathcal U}(N), {\mathcal L})$ is a
nondegenerate singularity of codimension
$k$ and Williamson type $(k_e,k_h,k_f)$,
of an IHS with $n$ degrees of freedom, then $k = k_e + k_h + 2 k_f$,
the ellipticity of $N$ is $k_e$, the hyperbolicity of $N$ is $k_h + k_f$,
and the
closedness of $N$ is $c(N) = n - k + k_f = n - k_e - k_h - k_f$. Orbits
of maximal dimension in $N$ are of dimension $n - k_e$, diffeomorphic to
${\bf T}^{c(N)} \times {\bf R}^{k_h + k_f}$ and
have orbit type $(k_e,0,0,c(N),k_h + k_f)$. It is also clear from the local
structure that the union of orbits of maximal dimension is dense in $N$.

\begin{thm}
\label{thm:action}
Let $({\mathcal U}(N), {\mathcal L})$ be a nondegenerate singularity of
Williamson type $(k_e,k_h,k_f)$ and codimension $k$,
of an IHS with $n$ degrees of freedom. Then there is
a natural Hamiltonian torus ${\bf T}^{c(N) + k_e}$ action in
$({\mathcal U}(N), {\mathcal L})$ which preserves the moment map of the IHS.
This action is unique, up to automorphisms of ${\bf T}^{c(N) + k_e}$, and it
is free almost everywhere in ${\mathcal U}(N)$. The isotropy group of this
action at $N$ (i.e. the subgroup of ${\bf T}^{c(N) + k_e}$, consisting of
elements which act trivially in $N$) is a torus ${\bf T}^{k_e}$ subgroup of
${\bf T}^{c(N) + k_e}$.
\end{thm}

{\it Proof}.
The proof uses the same arguments as in Theorem \ref{thm:Tn-1action} and
Theorem \ref{thm:focus}. Let ${\mathcal O}(x)$ be an orbit
maximal dimension in $N$.
Then ${\mathcal O}(x)$  has orbit type $(k_e,0,0,c(N),k_h + k_f)$. In
particular, it consists elliptic singular points (if $k_e > 0$). The
Poisson ${\bf R}^n$ action of our IHS on ${\mathcal O}(x)$ has a unique
${\bf R}^{n - k_h - k_f} = {\bf R}^{c(N) + k_e}$ subaction with compact
orbits. Denote the corresponding subgroup of ${\bf R}^n$ by
${\bf R}^{c(N) + k_e}_x$. By considering the local structure of smaller
orbits in the closure of ${\mathcal O}(x)$, we can see easily that the Poisson
${\bf R}^{c(N) + k_e}_x$ subaction has compact orbits on all $({\bf R}^n-)$
orbits of maximal dimension in $N$. Hence all orbits of this
${\bf R}^{c(N) + k_e}_x$ subaction in $N$ are compact (they are compact
subspaces of orbits of our IHS). We can redenote ${\bf R}^{c(N) + k_e}_x$
by ${\bf R}^{c(N) + k_e}_N$

If ${N'}$ is a regular leaf in ${\mathcal U}(N)$, then ${N'}$ passes
nearby at least one of the orbits of maximal dimension in $N$. From the
structure of elliptic singularities, it follows that there is a unique
$(c(N) + k_e)$-dimensional subspace ${\bf R}^{c(N) + k_e}_{N'}$ of
${\bf R}^n$, which is near to ${\bf R}^{c(N) + k_e}_x$, such that the
Poisson ${\bf R}^{c(N) + k_e}_{N'}$ subaction has compact orbits in
$\mathcal O$ (and is locally free there).

By continuity, we can associate to
each leaf ${N'}$ in ${\mathcal U}(N)$ (regular or not)  a subspace
${\bf R}^{c(N) + k_e}_{N'}$ of ${\bf R}^n$, in a continuous way, so that
the Poisson ${\bf R}^{c(N) + k_e}_{N'}$ subaction has compact orbits in
${N'}$, and is locally free if ${N'}$ is regular. Using
Arnold-Liouville  theorem, we can derive from this family of subspaces a
unique ${\bf T}^{c(N) + e_k}$ symplectic action which preserves the moment
map. The isotropy group of this action at $N$ corresponds to the
$k_e$-dimensional subspace of ${\bf R}^{c(N) + k_e}_x$ which acts trivially
on ${\mathcal O}_x$ (under the Poisson ${\bf R}^n$ action of our IHS).
$\blacksquare$

By the way, observe that the isotropic group at $N$ discussed above
arises from $k_e$ degrees of ellipticity of $N$, and this group acts ``around''
$N$ with general orbits being small-sized $k$-tori.

For minimal orbits in $N$ (i.e. closed orbits of dimension $n-k$), the
associated isotropy groups of the ${\bf T}^{c(N) + k_e}$ action are closed
subgroups of dimension $k_e + k_f$, containing the isotropy group at $N$.
If we chose a subgroup ${\bf T}^{n-k}_0$ of ${\bf T}^{c(N) + k_e}$, which is
transversal to all of the above isotropy groups, then the symplectic action
of ${\bf T}^{n-k}_0$ (as a subaction of ${\bf T}^{c(N) + k_e}$) is locally
free in ${\mathcal U}(N)$ and we can use Marsden-Weinstein reduction with
respect to this group action. But in order to avoid orbifolds, we need this
${\bf T}^{n-k}_0$ action to be free. In particular, ${\bf T}^{n-k}_0$ should
intersect each of the above $(k_e + k_f)$-dimensional isotropy groups only
at 0.

In order to achieve the freeness of ${\bf T}^{n-k}_0$ action, we
will use some covering argument like in Proposition \ref{prop:doublecover}
and some additional condition about isotropy groups like in Proposition
\ref{prop:stable2f}. Recall that a {\sc normal covering} of a topological
space $X$ is a covering $\overline{X}$
of $X$ which corresponds to a normal subgroup of
the fundamental group of $X$. In this case there is a free discrete group
action in $\overline{X}$ with the quotient space equal to $X$.

\begin{thm}
\label{thm:anycover} \footnote{erratum: Theorem \ref{thm:anycover} and its proof has
a hole: A-priori, there may be relations among $\alpha_i$, $\beta_i$ and $\gamma_j$
that have been omitted in the proof. For example, $\alpha_i$ may be homotopic to 0.
One way to show that such ``extra'' relations don't happen is to use the ``parallel
transport argument'' similar to the one used in the proof of Theorem
\ref{thm:aa-coordinates}. But this ``parallel transport argument'' requires the
topological stability condition or the absence of focus-focus components. Thus in
Theorem \ref{thm:anycover} we must require the singularity to be topologically
stable or not to contain focus-focus components. The finite covering is in general
not unique (situation similar to that of Proposition \ref{prop:doublecover}).} Let
${\mathcal U}(N)$ be as in theorem \ref{thm:action}. Then if the zero components of
the isotropy groups at minimal orbits of $N$ coincide, then there is a natural
finite normal covering $\overline{{\mathcal U}(N)}_{can}$ of ${\mathcal U}(N)$, such
that everything (symplectic form, moment map, ${\bf T}^{c(N) + k_e}$ action) can be
lifted from ${\mathcal U}(N)$ to $\overline{{\mathcal U}(N)}_{can}$, and the
isotropy groups at minimal orbits of the preimage $\overline{N}$ of $N$ are
connected and coincide.
\end{thm}

{\it Proof}. The proof is exactly the same as that of Theorem
\ref{thm:AL-1h}, and is based on the following presentation of the
fundamental group $\pi_1(N) = \pi_1 ({\mathcal U}(N))$ of $N$:

{\it Generators} (three types): $\alpha_1,...,\alpha_{n-k}$, which are
generated by an orbit of the ${\bf T}^{c(N) + k_e}$ action. (Notice that
there are only $n-k$ nonzero generating cycles of this type,
because the other $k_e + k_f$ cycles are vanishing due to focus-focus and
elliptic components). $\beta_1,...,\beta_s (s \geq 0)$,
which are ``exceptional cycles'': they lie in closed orbits of $N$ and are not
conjugate to an integral combination of $\alpha_i$. $\gamma_1,...,\gamma_t
(t \geq 1)$, which are ``base cycles'': the fundamental group
of the quotient space of $N$ by the ${\bf T}^{c(N) + k_e}$ action is
generated by the image of these cycles.

{\it Relations}: $\alpha_i$ commute with all the other generators. Twice of
$\beta_i$ are integral combinations of $\alpha_i$. If say, $\beta_1$ and
$\beta_2$ lie in the same closed orbit of $N$, then they commute. It may be
that some of $\beta_i$ which lie in different minimal orbits are conjugate
(if, say, they are conjugate to a same cycle which lies in an non-closed
orbit of $N$).

There is a natural homomorphism from $\pi_1(N)$ to ${\bf T}^{n-k}$ which maps
$\alpha_1,..., \alpha_{n-1},$ \ $\gamma_1,..., \gamma_t$ to 0 and $\beta_1,...,
\beta_s$ to elements of order 2. Take $G_{can}$ to be the kernel of this
homomorphism, and $\overline{{\mathcal U}(N)}_{can}$ to be the normal finite
covering of ${\mathcal U}(N)$ corresponding to $G$.

Like in Section \ref{section:focus}, the condition about isotropy groups
in Theorem \ref{thm:anycover} is
a kind of topological stability condition. The following general definition
of topological stability coincides with the ones given in the previous
sections for that particular cases discussed there.

\begin{define}
\label{define:stable}
A nondegenerate singularity $({\mathcal U}(N), {\mathcal L})$  of an IHS
is called {\sc topologically stable} if the local singular
value set of the moment map restricted to ${\mathcal U}(N)$ coincides with the
singular value set of the moment map restricted to a small neighborhood of
a singular point of maximal corank in $N$.
\end{define}

Because of Vey-Eliasson theorem about the local structure of singular
points, it is easy to describe the above local singular value sets.
Examples: 1) $x$ is a singular fixed point of Williamson type $(0,n,0)$, then
the local singular value set at $x$ of the moment map is a subset in
${\bf R}^n(F_1,...,F_n)$ diffeomorphic to the union of $n$
hyperplanes $\{F_i = 0 \}$. 2) If $x$ is a fixed point of Williamson type
$(0,1,1)$, then this set is diffeomorphic to the union of the line
$\{ F_2 = F_3 = 0 \}$ with the plane $\{ F_1 = 0 \}$ in the 3-space. 3)
If $x$ is an codimension 1 elliptic point then this set is diffeomorphic to
a closed half-space.

It is clear from the definition that elliptic singularities are
automatically topologically stable. It is as well
easy to construct examples of
nondegenerate hyperbolic and focus-focus singularities, that are not
topologically stable. In general, we have the following result, whose proof
is straightforward:

\begin{prop}
\label{prop:stable}
If $({\mathcal U}(N), {\mathcal L})$ is a nondegenerate topologically stable
singularity of codimension $k$ of an IHS with $n$ degrees of freedom
then we have: \\
a) All singular leaves of codimension $k$ in ${\mathcal U}(N)$ are topologically
equivalent. \\
b) A (germ of a) tubular neighborhood of any singular leaf in ${\mathcal U}(N)$
is a nondegenerate topologically stable singularity. \\
c) All closed orbits in $N$ have the same (minimal) dimension $n-k$. \\
\end{prop}

We mention the following important consequence of Theorem \ref{thm:anycover}
and Proposition \ref{prop:stable}:

\begin{cor}
If $({\mathcal U}(N), {\mathcal L})$ is a nondegenerate singularity of codimension
$k$ of an IHS with n degrees of freedom, which is topologically stable or
which has $k_f(N) =0$, then in an appropriate finite covering
$\overline{{\mathcal U}(N)}$ there is a free symplectic torus ${\bf T}^{n-k}$
action which preserves the moment map.
\end{cor}

\section{Topological decomposition}
\label{section:decomposition}
In this section all singularities are assumed to be topologically stable
nondegenerate. If $N_1$ and $N_2$ are two nondegenerate singular leaves in
two different IHS's, of codimension $k_1$ and $k_2$,
with the corresponding
Lagrangian foliation $({\mathcal U}(N_1), {\mathcal L}_1)$ and $({\mathcal U}(N_2), {\mathcal
L}_2)$, then the {\sc direct product} of these singularities
is the singular leaf $N = N_1 \times N_2$ of codimension $k_1 + k_2$ with
the associated Lagrangian foliation equal to the direct product of
the given Lagrangian foliations:
$$ ({\mathcal U}(N), {\mathcal L}) =
({\mathcal U}(N_1), {\mathcal L}_1) \times ({\mathcal U}(N_2), {\mathcal L}_2).$$

\begin{define}
A nondegenerate singularity $({\mathcal U}(N), {\mathcal L})$ of codimension $k$ and
Williamson type $(k_e, k_h, k_f)$ of an IHS with $n$ degrees of freedom
is called of {\sc direct-product type} topologically
(or a {\sc direct-product singularity}) if
it is homeomorphic, together with the Lagrangian foliation,
to a following direct product : \\
$ ({\mathcal U}(N), {\mathcal L}) \stackrel{homeo}{=} \\
({\mathcal U}({\bf T}^{n-k}), {\mathcal L}_r) \times
(P^2(N_{1}), {\mathcal L}_{1}) \times \ldots \times
(P^2(N_{k_e + k_h}), {\mathcal L}_{k_e + k_h}) \times
(P^4(N'_{1}), {\mathcal L}'_{1}) \times \ldots \times
(P^4(N'_{k_f}), {\mathcal L}'_{k_f}),$  \\
where
$({\mathcal U}({\bf T}^{n-k}), {\mathcal L}_r)$ denotes the Lagrangian foliation
in a tubular neighborhood of a regular Lagrangian $(n-k)$-torus of an IHS
with $n-k$ degrees of freedom;
$(P^2(N_{i}), {\mathcal L}_{i})$ for $1 \leq i \leq k_e + k_h$
denotes a codimension 1 nondegenerate
surface singularity (= singularity of an IHS with one degree of freedom);
$(P^4(N'_{i}), {\mathcal L}'_{i})$ for $1 \leq i \leq k_f$ denotes
a focus-focus singularity of an IHS with two degrees
of freedom.
\end{define}

\begin{define}
\label{define:adp}
A nondegenerate singularity of an IHS
is called of {\sc almost-direct-product type}
topologically (or simply an {\sc almost-direct-product singularity})
if a finite covering of it is homeomorphic, together with the Lagrangian
foliation, to a direct-product singularity.
\end{define}

Our main result is that any topologically
stable nondegenerate singularity is of almost direct product type. More
precisely, we have:

\begin{thm}
\label{thm:adp}
If $({\mathcal U}(N), {\mathcal L})$ is a nondegenerate topologically stable
singularity of William- \\
son type $(k_e, k_h, k_f)$ and codimension $k$ of an
IHS with $n$ degrees of freedom then it can be written
homeomorphically in the form of a quotient of a direct product
singularity \\
$({\mathcal U}({\bf T}^{n-k}), {\mathcal L}_r) \times
(P^2(N_{1}), {\mathcal L}_{1}) \times \ldots \times
(P^2(N_{k_e + k_h}), {\mathcal L}_{k_e + k_h}) \times
(P^4(N'_{1}), {\mathcal L}'_{1}) \times \ldots \times
(P^4(N'_{k_f}), {\mathcal L}'_{k_f}),$  \\
by a free action of a finite group $\Gamma$ with the following property:
$\Gamma$ acts on the above product
component-wise (i.e. it commutes with the projections onto the components),
and moreover, it acts trivially on elliptic components.
\end{thm}

{\bf Remark}. The above decomposition is not symplectic, i.e. in general
we cannot decompose the symplectic form to a direct sum of
symplectic forms of the components. However, one will see from the proof
that we can replace the word
``homeomorphically'' by the word ``diffeomorphically''.

A direct product with an action group in the above theorem will be called a
{\sc model} of a stable nondegenerate singularity
$({\mathcal U}(N), {\mathcal L})$. A model is called {\sc canonical} if there does not
exist a nontrivial element of $\Gamma$ which acts trivially on all of the
components except one.

\begin{prop}
\label{prop:minimalmodel} If $({\mathcal U}(N), {\mathcal L})$ is a topologically
stable nondegenerate singularity  then there is a unique\footnote{erratum: not
unique in general, but unique if $n=k$, i.e. if $N$ contains a fixed point}
canonical model $(\overline{{\mathcal U}(N)}_{CAN}, {\rm action \; of} \;
\Gamma_{CAN})$ for it.
\end{prop}

The group $\Gamma_{CAN}$ which enters in the canonical model will be called the {\sc
Galois group} of the singularity $({\mathcal U}(N), {\mathcal L})$. remark that for
codimension 1 singularities, $\overline{{\mathcal U}(N)}_{CAN} = \overline{{\mathcal
U}(N)}_{can}$. In general, $\overline{{\mathcal U}(N)}_{CAN}$ is a finite covering
of $\overline{{\mathcal U}(N)}_{can}$\footnote{Remember that $\overline{{\mathcal
U}(N)}_{can}$ is not unique in general, though we can still keep that notation. It
seems that if $k <n$ then in general the group $\Gamma_{CAN}$ is not unique either.
}.

\begin{cor}
i)\footnote{erratum: in i), must add the condition that $N$ does not contain
focus-focus components} If $N$ is a topologically stable nondegenerate singular leaf
then it is an
Eilenberg-Maclane space $K(\pi_1(N), 1)$. \\
ii)\footnote{erratum: ii) is not proved and probably not true} If $\phi: N \to N$ is
a homeomorphism from $N$ to itself whose induced automorphism on $H^1(N, {\bf R}$ is
identity, then
$\phi$ is isomorphic to identity. \\
iii)\footnote{addendum: perhaps one can use ``collapsing arguments'' to extract
information about $N$ from the boundary of its tubular neighborhood and to prove
iii)} If two topologically stable nondegenerate singularities of codimension greater
than 1 are such that their nearby singularities of smaller codimensions are
topologically equivalent in a natural way, then these two singularities are also
topologically equivalent.
\end{cor}

The above corollary answers a question posed in \cite{Zung94} and is usefull
in the computation of the canonical model of singularities of well-known
integrable systems.

{\it Proof of Theorem \ref{thm:adp} and Proposition \ref{prop:minimalmodel}}. We
will prove for the case $k=n$ (i.e. the case when $N$ contains a fixed point). The
case $k < n$ can be then easily proved by the use of the results in Section
\ref{section:action} and Section \ref{section:aa-coordinates}.( In case $k < n$, we
will use finite covering twice, but the construction is canonical\footnote{erratum:
the construction is not unique when $k < n$, though there are only a finite number
of possibilities}, and one can verify directly from the construction of subgroups of
$\pi_1(N)$ that we still have a normal covering).

Let $({\mathcal U}(N), {\mathcal L})$ have maximal possible codimension $k = n$.
For simplicity, we will also assume that $N$ has no elliptic component, i.e.
$k_e(N) = 0$ (the case where $k_e(N) > 0$ can then be proved with the aid of
the torus action around $N$ which arises from elliptic components (see
Section \ref{section:action})). We will use the following notion of
$l-type$ of $({\mathcal U}(N), {\mathcal L})$,
which was introduced by Bolsinov for the case of hyperbolic
singularities of IHS with two degrees of freedom.

Let $x$ be a fixed point in $N$.  By changing the Poisson action but leaving
the topological structure of $({\mathcal U}(N), {\mathcal L})$ unchanged, we can
assume that at point $x$ the moment map is linear and its components are as
in Williamson's Theorem \ref{thm:Williamson}. We will denote the hyperbolic
components of the moment map by $F_1,...,F_{k_h}$, its focus-focus
components by $F'_1,F''_2,...,F'_{k_f},F''_{k_f}$ $(k_h + 2 k_f = n)$. Then
we have $F_i = x_i y_i, F'_i = x'_i y''_i - x''_i y'_i, F''_i = x'_i y'_i +
x''_i y''_i$, where $(x_1,y_1,...,x_{k_h}, y_{k_h}, x'_1, y'_1, x''_1,
y''_1,...,x'_{k_f}, y'_{k_f}, x''_{k_f}, y''_{k_f})$ is a canonical system
of coordinates at $x$. The local singular value set of our moment map
restricted to a small neighborhood of $x$ will be a germ at zero of the
union of codimension 1 hyperplanes $\{ F_i = 0 \}$ and codimension 2
hyperplanes $\{ F'_i = F''_i = 0 \}$ in ${\bf R}^n$.
By definition of topological stability, this set is also the singular value
set of the moment map restricted to ${\mathcal U}(N)$. Denote by $I_i, 1\leq i
\leq k_h$ (resp., $k_h < i \leq n$)
the subset of this singular value set, which have all
coordinates equal to zero except $F_i$ (resp., $F'_i$ and $F''_i$).
Let $V_i$ $(1 \leq i \leq n)$ denotes the preimage of $I_i$ in ${\mathcal U}(N)$
of the projection. Then $V_i$ are symplectic 2-manifolds ($1 \leq i \leq
k_h$) and symplectic 4-manifolds ($k_h < i \leq n$),
and the intersection of the Lagrangian foliation
in ${\mathcal U}(N)$ with $V_i$ gives a singular Lagrangian foliation
to these manifolds. With this Lagrangian foliation structure,
$V_i$  become
hyperbolic and focus-focus singularities of IHS's with one and two degrees
of freedom, respectively. We will denote the
singularities associated to $V_i$ by $V_i$ again.
Strictly speaking, since $V_i$ may be non-connected, it may be
that they are not singularities but a finite set of singularities.
The (unordered) $n$-tuple $(V_1,...,V_n)$
will be called the {\sc l-type} of the singularity $({\mathcal U}(N), {\mathcal
L})$. (In case $k_e > 0$, we add elliptic codimension 1 singularities of IHS
with one degree of freedom to this l-type).

Denote the singular leaf (or more precisely, the union of singular leaves)
in $V_i$ by $K_i$. Then $K_i$ belong to $N$.
Recall that $N$ has a natural stratification, given by its orbits. The union
of $K_i$, $1 \leq i \leq k_h$, is the 1-skeleton of $N$.
The union of all $K_i$, which we will denote by $Spine(N)$, is the
two-skeleton of $N$ minus open 2-dimensional orbits (orbits diffeomorphic to
${\bf R}^2$). Observe the following important fact: $Spine(N)$ gives all
generators of the fundamental group $\pi_1(N)$. Other orbits in $N$ don't
give any new cycles of $\pi_1(N)$, just new commutation relations.

By a {\sc primitive orbit} in $N$ we will mean an
one-dimensional orbit, or a two-dimensional orbit of
one degree of closedness and one degree of openness. In other words,
a primitive orbit is an orbit that lies in $Spine(N)$ and is not a
fixed point. Let ${\mathcal O}_i (i = 1,2)$ be a primitive orbit
with two boundary (limit) points $x_i$ and $y_i$, such that $x_2 = y_1$.
For each $i$, $x_i$ and $y_i$ are fixed points in $N$, which
may coincide. Provide ${\mathcal O}_i$ with the direction going from $x_i$ to
$y_i$ so that they become oriented 1-cells.
Suppose that ${\mathcal O}_1$ and ${\mathcal O}_2$ lie in different
$V_i$ (components of l-type). Then we will construct a mapping from
the oriented orbit ${\mathcal O}_1$ to another oriented orbit ${\mathcal O}'_1$
which belongs to the same $V_i$
${\mathcal O}_1$. This mapping will be well defined, and
we will call it the
{\sc moving of ${\mathcal O}_1$ along ${\mathcal O}_2$}.

We will define the above moving in the case when ${\mathcal O}_i$ are
one-dimensional (i.e. they lie in $V_i$ but not $V'_i$).
The other cases are similar. (To see why, take the
quotient of ${\mathcal U}(N)$ by the ${\bf T}^{k_f}$ action arisen from
focus-focus components as in Section \ref{section:action}). It follows from
the local structure theorems that there is an open 2-dimensional
orbit (diffeomorphic to ${\bf R}^2$) in $N$, denoted by $\mathcal O$, which
contains ${\mathcal O}_1$ and ${\mathcal O}_2$ in its closure. The algebraic boundary
of $\mathcal O$ consists of four 1-cells (1-dimensional orbits)
and  four 0-cells. Denote the other two 1-cells by ${\mathcal O}'_1$,
${\mathcal O}'_2$, in such an order that $\partial {\mathcal O} = {\mathcal O}_1 +
{\mathcal O}_2 - {\mathcal O}'_1 - {\mathcal O}'_2$ algebraically.
Then it is easy to be seen (using the moment map)
that ${\mathcal O}'_1$ belongs to the same $V_i$ or
$V'_i$ as ${\mathcal O}_1$. The moving of ${\mathcal O}_1$
along ${\mathcal O}_2$ is by definition the oriented orbit
${\mathcal O}'_1$ (with the end point at $y_2$). This moving can also be
understood as an orientation homeomorphism from ${\mathcal O}_1$ to
${\mathcal O}'_1$, which maps end points to end points, and which is defined up
to isotopies. The pair oriented $({\mathcal O}'_2,{\mathcal O}'_1)$ will
be called the {\sc elementary homotopy deformation} of the oriented
pair $({\mathcal O}_1, {\mathcal O}_2)$.

If ${\mathcal O}_i$ (i = 1,...,s) is a chain of oriented primitive orbits in $N$
with start and end points $x_i, y_i$ such that $x_i = y_{i-1}, x_1 = y_s$,
then they form a closed 1-cycle in $N$. (Here, say if ${\mathcal O}_1$ lies in
some focus-focus $V_i$, they we replace ${\mathcal O}_1$ by a curve lying in
${\mathcal O}_1$ and going from $x_1$ to $y_1$).  We know that any closed curve
in $N$ is homotopic to such a cycle. Furthermore, any two homotopic cycles
can be obtained from one another by a finite number of elementary homotopy
deformations. That is because orbits with three or more degrees of openness
in $N$ do not contribute any new generator or relation to $\pi_1(N)$.

From now on, in our homotopy arguments by a closed curve
in $N$ we will always mean  a cycle consisting of consecutive primitive
orbit. Let $x$ be a fixed point in $N$. Let $W_i$, $1 \leq i \leq n$,
be the connected component of $V_i$  that contains $x$.
Let ${\gamma}$ be a closed curve
going from $x$ to $x$ and
lying entirely in $V_1$. Let ${\mathcal O}$ be an oriented primitive orbit with
the start point at $x$ and not lying in $V_1$.
Then we can move ${\mathcal O}$ along
$\gamma$ by moving it step by step along orbits in $\gamma$. It is clear
that the end result will depend only on the homotopy type of $\gamma$ as an
element in $\pi_1(W_1, x)$. Denote this end result by ${\mathcal O}^{\gamma}$.
It is equally clear that $x$ is also the start point of ${\mathcal O}^{\gamma}$.
If ${\mathcal O}$ belongs to $W_i$ ($i \neq 1$), say, then there are 4 local
primitive orbits lying in $W_i$ having $x$ as the start point. One can move
all of these local orbits along $\gamma$. In other words, one can move a
local singular leaf of $x$ in $W_i$ along $\gamma$. Recall that $W_i$ has a
natural orientation given by the induced symplectic structure. In can be
easily seen by induction along primitive orbits in $\gamma$ that the above
moving is orientation-preserving local homeomorphism of the local singular
leaf of $x$ in $W_i$, which is isotopic either to identity or to an
involution that maps each local primitive orbit in $W_i$ at $x$ to its
opposite. It follows that the subset of $\pi_1(W_1)$ consisting of (homotopy
type of) closed curves whose action on each of the local singular leaf of
$x$ in $W_i$ is identity (in isotopy category) for all $i \neq 1$, is a
normal subgroup of $\pi_1(W_1)$ of index at most $2^{n-1}$. Denote the image
of this subgroup in $\pi_1(N, x)$ by $G_1$.

By replacing $W_1$ by other $W_i$, we can construct subgroups $G_i$ ($1 \leq
i \leq n$) of $\pi_1(N, x)$ in the same way. A very important property of
$G_i$ that we will show below can be stated roughly as follows: the moving
along elements of $G_1$ (1 can be replaced by any number from 1 to $n$)
gives the trivial action of $G_1$ not only on the local singular leaf of $x$
in $W_i$ ($i \neq 1$),  but also on the whole singular leaf through $x$ in
$W_i$, i.e. the intersection of $W_i$ and $K_i$.

Let ${\mathcal O}_1,...,{\mathcal O}_s$ be a chain of consecutive primitive orbits
lying in $W_i$, $i \neq 1$, and $x_1 = x, y_1 = x_2,..., y_s$ their
corresponding start and end points. Let $\gamma$ be a closed curve in $W_1$
which represents an element in $G_1$. Then $\gamma$ moves ${\mathcal O}_1$ to
itself. We can also move $\gamma$, considered not as a closed curve but
simply a curve with  the start point at $x$, along ${\mathcal O}_1$. Let the
image be called $\gamma_1$. Then $\gamma_1$ will be in fact a closed curve
starting and ending at $x_2$ (exactly because $\gamma$ moves ${\mathcal O}_1$ to
itself). Moreover, $\gamma_1$ moves $- {\mathcal O}_1$ (considered with the
inverse direction) to itself. As a consequence, $\gamma_1$ also moves ${\mathcal
O}_2$ to itself, and we can continue like this. In the end, we have that the
curve $({\mathcal O}_1,...,{\mathcal O}_s)$ is moved to itself along $\gamma$. Here
we should note that the motion of $({\mathcal O}_1,...,{\mathcal O}_s)$ along any
curve $\gamma$ (closed or not) with the starting point at $x$ and with no
component belonging to $W_i$ is always well defined by applying step by step
the elementary homotopy deformations and using the definition of moving for
primitive orbits.

The main consequence of the above trivial action property of $G_i$ is as
follows: if $\alpha \in G_i$ and $\beta \in G_j$ with $i \neq j$ then
$\alpha$ and $\beta$ commute. In other words, the subgroups $G_i$ of
$\pi_1(N)$ commute pairwise. Indeed, let $\alpha$ be a closed curve in $W_1$
representing an element in $G_1$, $\beta$ a closed curve in $W_2$
representing an element in $G_2$. Then $\alpha$ moves $\beta$ to a curve
${\beta}'$ which is in fact equal to $\beta$, and $\beta$ moves
$\alpha$ to a curve ${\alpha}'$ which is in fact equal to $\alpha$, and we
have $\alpha \beta = {\beta}' {\alpha}'$ (because by definition,
${\beta}' {\alpha}'$ is obtained from ${\beta}{\alpha}$ by a finite number
of elementary homotopy deformation).

Let $G$ denote the product of $G_i$ in $\pi_1(N)$.
We will now show that $G$ is a normal subgroup of finite index in
$\pi_1(N)$.

First, about the normality of $G$. We will call it $G(x)$, to remember that it
is a subgroup of $\pi_1(N, x)$. If $y$ is another fixed point of $N$, then
similarly we can construct the subgroup $G(y)$ of $\pi_1(N, y)$. Let
$\beta = ({\mathcal O}_1,...,{\mathcal O}_s)$
be a chain of consecutive primitive orbits
starting at $x$ and ending at $y$, we will construct a natural homomorphism
from $G(x)$ to $G(y)$, called
{\sc moving from $G(x)$ to $G(y)$ along $\beta$}.
By going back from $y$ to $x$ via another chain ${\beta}'$, we will get a
similar homomorphism from $G(y)$ to $G(x)$. Combining these two homomorphisms,
we will get an automorphism of $G(x)$, which is obtained in fact by the
conjugation of $G(x)$ with the cycle $\beta {\beta}'$. Thus the conjugation
of $G(x)$ with any cycle in $\pi_1(N, x)$ is $G(x)$ itself, i.e. $G(x)$ is
a normal subgroup of $\pi_1(N, x)$.

To construct the above moving from $G(x)$ to $G(y)$, we can assume for
simplicity that the chosen chain $\beta$ from $x$ to $y$ consists of only one
orbit ${\mathcal O}_1$. ${\mathcal O}_1$ lies in some $W_i$, say $W_1$ for
definiteness. Then if $\gamma \in G(x)$ is represented by a closed curve in
$W_i, i \neq 1$, it is moved to an element in $\pi_1(N, y)$ by moving of the
curve along ${\mathcal O}_1$ in the usual way as before. If $\gamma$ is
represented by a closed curve in $W_1$, then it is moved to the element of
$\pi_1(N)$ which is represented by the closed curve $(- {\mathcal O}_1 +
\gamma + {\mathcal O}_1)$. If $\gamma$ is a product of the above generators of
$G(x)$, then it is moved to the product of the associated images. Because of
the trivial action property of  $G(x)$, one can easily check that the above
moving is well defined and will be a homomorphism from $G(x)$ to $G(y)$.

Now, about the finiteness of the index of $G$ in $\pi_1(N)$. Let $\gamma$ be
an arbitrary element of $\pi_1(N, x)$, which is represented by a closed curve
(consisting of primitive orbits). Using a finite number of elementary
homotopy deformations, we see that $\gamma$ can be represented as
$\beta_1 + ... + \beta_n$, where each $\beta_i$ is a chain of consecutive
primitive orbits lying entirely in $W_i$. Because of the number of
primitive orbits in $N$ is finite, and the index of $G_i$ in $W_i$
(for each fixed $x$) is finite, we can chose for each $i$ a finite number of
curves (i.e. chains of consecutive primitive orbits)
$\beta^1_i,..., \beta^i_I$, such that any curve $\beta_i$ in $V_i$ is
homotopic (rel. end points) to $\gamma_i \beta_i^{s(i)}$, where $\gamma_i$
is an element in some group $G(z)$ (i.e. freely conjugate to an element of
$G(x)$), and $1 \leq s(i) \leq I$. Then
$\gamma$ is homotopic to $\gamma_1 + \beta_1^s(1) + ... +
\gamma_n + \beta_n^{s(n)}$. Now, because $\gamma_n$ has the trivial action
property, $\beta_{n-1}^{s(n-1)} + \gamma_n$ is homotopic (rel. end points) to
$\gamma_n' + \beta_{n-1}^{s(n-1)}$, where $\gamma_n'$ again belongs to some
$G(z)$. By induction, we obtain that $\gamma$ is homotopic to
$\gamma'_1 + \beta_1^{s(1)} + ... + \beta_n^{s(n)}$, where $\gamma_1'$ is some
element of $G(x)$. This proves the finiteness of the index of $G$ in
$\pi_1(N)$.

Let $\overline{{\mathcal U}(n)}$ denotes the normal finite covering of
${\mathcal U}(N)$ associated to the above normal subgroup $G$ of the
fundamental group   $\pi_1({\mathcal U}(N)) = \pi_1 (N)$. It is obvious that
the symplectic structure, the singular Lagrangian foliation and the moment map
can be lifted from ${{\mathcal U}(n)}$ to $\overline{{\mathcal U}(n)}$.
$\overline{{\mathcal U}(n)}$ has the fundamental group equal to $G$, and this group
has the trivial action property inherited from $N$. From this one can easily
show that $\overline{{\mathcal U}(n)}$ has the direct-product type.
Let $\Gamma$ be the quotient group $\pi_1(N,x) / G(x)$. Then $\Gamma$ acts
freely on $\overline{{\mathcal U}(N)}$, with the quotient being the singularity
$({\mathcal U}(N), {\mathcal L})$. This action can be made component-wise, because
it $\Gamma$ acts on the associated finite covering of $Spine(N)$
component-wisely. (It is an easy exercise to see why).
Theorem \ref{thm:adp} is proved.

It is also clear that only elements of $G$ have the trivial action property.
As a consequence, if $G'$ is another subgroup of $\pi_1(N, x)$ such that
the associated finite covering is of product type, then $G'$ must be a
subgroup of $G$. It follows that the model constructed above (for $k=n$)
is the unique canonical model. Proposition \ref{prop:minimalmodel} is proved.
$\blacksquare$

{\bf Example}s. 1) Lerman and Umanskii \cite{LU2} (see also \cite{Bolsinov})
classified topologically stable nondegenerate
hyperbolic codimension 2 singularities of IHS's with two degrees with freedom,
which contain only one fixed point. The result is that, topologically there
are four different cases. But their description is rather complicated, in
terms of cell decomposition. It is an easy exercise to write down explicitly
the canonical model for all of these four cases (it was done in
\cite{ZungThesis}). 2) Computation of the canonical model for singularities of
the geodesic flow on multi-dimensional ellipsoids was done in
\cite{ZungSphere}.

\section{Action-angle coordinates}
\label{section:aa-coordinates}
The following theorem gives a non-complete system of action-angle
coordinates to each topologically stable nondegenerate singularity.
We will consider only the case when there is no elliptic component.
In case when there are some elliptic components, there will be even more
action-angle coordinates (some of which are polar).
This case is left to the reader.

\begin{thm}
\label{thm:aa-coordinates}
Let $({\mathcal U}(L), {\mathcal L})$ be a topologically stable nondegenerate
singularity of codimension $k$ and Williamson type $(0, k_h, k_f)$ of an IHS
with $n$ degrees of freedom. Suppose that a moment map preserving hamiltonian
${\bf T}^{n-k}$ action is free in ${\mathcal U}(N)$. Then we have: \\
a) The above action provides ${\mathcal U}(N)$ with a principal ${\bf T}^{n-k}$
bundle structure, which is topologically trivial. \\
b) There is a coisotropic section to this trivial bundle. \\
c) ${\mathcal U}(N)$ is symplectomorphic to the direct product ${\bf D}^{n-k}
\times {\bf T}^{n-k} \times P^{2k}$ with the symplectic form
$$\omega = \sum^{n-k}_1 dx_i \wedge d y_i + \pi^*(\omega_1)$$
where $x_i$ are Euclidean coordinates on ${\bf D}^{n-k}$, $y_i$
($ {\rm mod} \; 1$) are coordinates on ${\bf T}^{n-k}$, $\omega_1$ is a
symplectic form on a $2k$-dimensional symplectic manifold $P^{2k}$, and
$\pi$ means the projection. Under this symplectomorphism, the moment map does
not depend on $y_i$. The set of
functions $x_i, y_i$ will be called a {\sc non-complete
system of action-angle coordinates} for singularity
$({\mathcal U}(L), {\mathcal L})$. \\
d) If $({\mathcal U}(N), {\mathcal L})$ is any topologically stable nondegenerate
singularity of codimension k, then $\overline{{\mathcal U}(N)}_{can}$ (as in Theorem
\ref{thm:anycover}) is symplectomorphic to the direct product ${\bf D}^{n-k} \times
{\bf T}^{n-k} \times P^{2k}$ with the canonical symplectic form as before, and
moreover $\Gamma_{can}$ acts on this product component-wise. In other words, we have
an {\sc equivariant non-complete system of action-angle coordinates} on
$\overline{{\mathcal U}(N)}_{can}$ for each topologically stable nondegenerate
singularity.
\end{thm}

{\it Proof}.
a) By changing the moment map but leaving the Lagrangian foliation
unchanged, we can assume that the components $F_i , i = 1,..., n$ of the
moment map are chosen so that the Hamiltonian vector fields
$X_{n-k+1} = X_{F_{n-k+1}},...,X_{n} = X_{F_{n}}$
generate a free ${\bf T}^{n-k}$
action in ${\mathcal U}(N)$, which give rise to the principal bundle in question,
and the local singular value set of the moment map is given as the union of
local codimension 1 hyperplanes $\{ F_i = 0 \}$, $1 \leq i \leq k_h$
and codimension 2 hyperplanes $\{ F_i = F_i + 1 = 0 \}$,
$k_h < i \leq k, (k - i - 1) \vdots 2$.

Because the classifying space of ${\bf T}^{n-k}$ is the product of $(n-k)$ samples
of ${\bf CP}^{\infty}$, a principal ${\bf T}^{n-k}$ bundle is trivial if and only if
a finite covering of it is also trivial. Thus, in view of Theorem \ref{thm:adp} (for
the case $n=k$), we can assume that the base space of ${\mathcal U}(N)$ is a union
of a $(n-k)$-dimensional family of topologically equivalent singularities of direct
product type. Namely, we will assume that the base space of ${\mathcal U}(N)$,
denoted by $B$, is homeomorphic to the direct product ${\bf D}^{n-k} \times W_1
\times ... \times W_{k_h + k_f}$, where ${\bf D}^{n-k}$ is coordinated by $F_1,...,
F_{n-k}$; each $W_i$ ($1 \leq i \leq k_h)$ is a hyperbolic singularity of an IHS
with one degree of freedom; and each $W_i$ ($k_h \leq i \leq k_h + k_f)$ is a
topologically stable focus-focus singularity of an IHS with two degrees of freedom.
Under this homeomorphism, $F_i$, $1 \leq i \leq k_h$, becomes a function depending
only on $W_i$, and $F_{k_h + 2i - 1},
 F_{k_h + 2i}$, $1 \leq i \leq k_f$, become functions depending only on
$W_{k_h + i}$.

Because of homotopy, it is enough to prove the triviality of the bundle in $N$.
Denote the base space of this bundle by $Spine(B)$ (it is the spine of $B$ in the
same sense as in Section \ref{section:decomposition}). We will define a parallel
transportation along the curves of $Spine(B)$ of the elements of $N$. This parallel
transportation will be locally flat (zero curvature), from which follows easily the
triviality of the ${\bf T}^{n-k}$ bundle structure of $N$. To define this
transportation, we simply use the vector fields $X_i = X_{F_i}$, $ 1 \leq i \leq k$.
First of all, notice that these vector fields commute. Second, where they are
different from zero, they are transversal to the tori of the bundle. Third, suppose
for example that $X_1 (x) = 0$ for some singular point $x \in N$. Then because of
our choice of the moment map (more precisely, because we choose the local singular
value set to be ``canonical''), the direction of $X_i$ near $x$ is transversal to
the orbit through $x$. The second property means that we have a parallel
transportation in regular orbits of $N$. The third property means that we can extend
this parallel transportation over the singular orbits of $N$, to obtain a well
defined parallel transportation (connection) in $N$. The first property means that
this connection is locally flat. Assertion a) is proved.

Assertion b) is a direct consequence of assertion c), because if one have a
non-complete system of action-angle coordinates then the submanifold
$\{ y_1 =...= y_{n-k} = 0 \}$ is obviously a coisotropic section to the
${\bf T}^{n-k}$ bundle. In fact, one can see easily that assertions b) and
c) are equivalent.

c) The proof of assertion c) is similar to that of \ref{thm:canonical1h}.
For completeness, we will recall it here.

Put $x_i = F_{i + k}, 1 \leq i \leq n-k$, where we suppose that the moment
map is chosen as in the proof of assertion a), and denote the Hamiltonian
vector field of $x_i$ by $\xi_i$..
Recall from assertion a) that the ${\bf T}^{n-k}$ bundle structure of
${\mathcal U}(N)$ is trivial. Let $L$ be an arbitrary section to this bundle
($L$ is diffeomorphic to the base space $B$),
and define functions $z_i \; (i = 1,\ldots,n-k)$
by putting them equal to zero on $L$ and setting
$1 = dz_i(\xi_i) = \{x_i, z_i\}$ (Poisson bracket). Set
$\omega_1 = \omega - \sum_1^{n-k} dx_i \wedge dz_i$.
Then one checks that $L_{\xi_i} \omega_1 = i_{\xi_i} \omega_1 = 0$. It means
that $\omega_1$ is a lift of some closed 2-form from $B$ to
${\mathcal U}(N)$, which we will also denote by $\omega_1$. Since $\omega$
is non-degenerate, it follows that $\omega_1$ is non-degenerate on every
$2k$-dimensional manifold (with boundary) $P^{2k}_{x_1,\ldots,x_{n-1}} =
B \cap \{x_1,\ldots,x_{n-1} {\tt fixed}\}$. Using Moser's path method
\cite{Moser}, one can construct a diffeomorphism
$\phi : B \to {\bf D}^{n-k} \times P^{2k}$, under which $\omega_1$
restricted to $P^{2k}_{x_1,\ldots,x_{n-k}}$ does not depend on the choice of
$x_1,\ldots,x_{n-1}$. In other words, there is a symplectic form $\omega_2$
on $P^{2k}$ such that $\omega_1 - \omega_2$ vanishes on every
$P^{2k}_{x_1,\ldots,x_{n-k}}$. Since $d(\omega_1 - \omega_2) = 0$, we can write
it as $\omega_1 - \omega_2 = d(\sum_1^{n-k}
a_i dx_i + \beta)$, where $\beta$ is
some 1-form on $B^{n+1}$ (which is not zero on $P^{2k}_{x_1,\ldots,x_{n-k}}$
in general). If we can eliminate
$\beta$, i.e. write $\omega_1 - \omega_2 = d(\sum_1^{n-k} a_i dx_i)$,
then we will have $\omega = (\sum_1^{n-k}
dx_i \wedge dz_i - \sum_1^{n-k} dx_i \wedge da_i + \omega_2
= \sum_1^{n-k}dx_i \wedge d(z_i - a_i) + \omega_2$, and the theorem will
be proved by  putting $y_i = z_i - a_i$. Let us show now how to
eliminate $\beta$. $\beta$ restricted on every $P^{2k}_{x_1,\ldots,x_{n-k}}$
is a closed 1-form, hence it represents a cohomology element
$[\beta](x_1,\ldots,x_{n-k}) \in H^1(P^{2k})$. If
$[\beta](x_1,\ldots,x_{n-k}) \equiv 0$ then
$\beta = dF - b_1 dx_1 - \ldots - b_{n-k} dx_{n-k}$
for some functions $F, b_1,\ldots,b_n$, and we have
$\omega_1 - \omega_2 = d(\sum(a_i - b_i x_i) dx_i)$. In general, we can
achieve $[\beta](x_1,\ldots,x_{n-k}) \equiv 0$ by induction on the number of
generators of $H^1(P^{2k}, {\bf R})$ as follows. Let $\gamma$ be
an element in a fixed system of generators of $H^1(P^{2k}, {\bf R})$. Set
$b(x_1,\ldots,x_{n-k}) = <[\beta], \gamma>(x_1,\ldots,x_{n-k})$.
Chose a closed 1-form $\delta$ on $P^{2k}$ such that $\delta (\gamma) = 1$
and the action of $\delta$ on other generators from the chosen system is 0.
Change $\omega_2$ for the following 2-form on ${\mathcal U}(N)$:
$\omega_2' = \omega_2 + db \wedge \delta$. It is clear that $\omega_2'$ and
$\omega_2$ restricted on every $P^{2k}_{x_1,\ldots,x_{n-k}}$ are the same.
Moreover, $\omega_2'$ is closed and of rank $2k$.
Thus the distribution by its
$(n-k)$-dimensional tangent zero-subspaces is integrable, and it gives rise
again to a diffeomorphism $\phi' : B \to {\bf D}^{n-k} \times P^{2k}$.
Replacing $\omega_2$ by $\omega_2'$, we have
$\omega_1 - \omega_2' = d(\sum_1^{n-k}
a_i dx_i + \beta')$, with $\beta' = \beta - b \delta$, whence
$<[\beta'],\gamma> = 0$.

d) The proof of d) is similar to that of Theorem \ref{thm:AL-1h}.
$\blacksquare$

\section{Nondegenerate IHS's}
\label{section:nondesys}
There has been no commonly accepted notion of nondegenerate IHS, though intuitively
one understands that such a notion should be connected to the nondegeneracy of
singularities. In this last Section we propose two extreme notions: strong
nondegeneracy and weak nondegeneracy for IHS's. We also suggest a middle condition,
which will be simply called nondegeneracy, and which seems to us most suitable for
practical purposes. Most, if not all, known IHS's in mechanics and physics satisfy
the (middle) condition of nondegeneracy. Some, but not all, known IHS's satisfy the
strong condition of nondegeneracy.

\begin{define}
An IHS is called {\sc strongly nondegenerate} if all of its
singularities are stable nondegenerate.
\end{define}

Examples of systems satisfying the strong nondegeneracy are the Euler and
Lagrange tops, and the geodesic flow on the multi-dimensional ellipsoid
where one cuts out the zero section in the cotangent bundle. Examples of
systems non satisfying strong nondegeneracy condition include the
Kovalevskaya top.

\begin{define}
An IHS is called {\sc weakly nondegenerate}
if the orbit space of its associated
singular Lagrangian foliation is Hausdorff, and almost all of its
singular points belong to nondegenerate singular leaves.
\end{define}

In the above definition, the word ``almost all'' means a dense subset in the
set of singular points. Since all singular leaves near a nondegenerate
singular leaf are also nondegenerate, we have in fact an open dense subset.

Let $x \in M^{2n}$ be a singular point of corank an IHS with the moment map
${\bf F} = (F_1,...,F_n)$. We can assume that $dF_1 \wedge ... \wedge
dF_{n-k} (x) \neq 0$. Applying a local Marsden-Weistein reduction with
respect to $F_1,...,F_{n-1}$, we obtain a local integrable system with $k$
degrees of freedom, for which $x$ becomes a fixed point. $x$ is called a
{\sc clean} singular point, if it becomes an isolated fixed point under this
local reduction. Obviously, all nondegenerate singular points are clean.

\begin{define}
An IHS is called {\sc nondegenerate} if it is weakly nondegenerate, all of
its nondegenerate singularities are topologically stable, and all of its
singular points are clean.
\end{define}

In a nondegenerate IHS, singular leaves which are not nondegenerate
will be called {\sc simply-degenerate}. Thus each singular leaf of a
nondegenerate IHS is either nondegenerate or simply-degenerate.
It can be shown easily by continuation that torus
actions discussed in Section \ref{section:action} still exist for
simply-degenerate singularities: in a saturated tubular neighborhood of a
simply-degenerate singularity of codimension $k$, there is a
symplectic $(n-k)$-dimensional torus action which preserves the moment map.
Moreover, an (equivariant) system of $(n-k)$ action and $(n-k)$ angle
coordinates still exist. It is perhaps the most important property of
simply-degenerate singularities.

One way to obtain (strongly or weakly) nondegenerate IHS's is to use {\sc
integrable surgery}, which will be discussed in more detail in
\cite{ZungSurgery}. Integrable surgery simply means the cutting and gluing
of pieces of symplectic manifolds together with IHS's on them so that to
obtain new symplectic structures and new IHS's.
In other words, one does a surgery on the level of $n$-dimensional orbit
spaces and tries to lift this surgery to the level of symplectic manifolds.
In principle, by this way one should be able to see if a stratified affine
manifold is an orbit space of any IHS, and how many IHS's does it correspond
to. Surprisingly or not, this
method gives a very simple way to construct many interesting symplectic
manifolds, including nonstandard symplectic ${\bf R}^{2n}$ and
Kodaira-Thurston example of a non-K\"ahlerian closed symplectic manifold.
For example, consider a simplest case, where the orbit space is a closed
annulus, whose boundary corresponds to elliptic singularities. Interestingly
enough, there is a discrete 2-dimensional family of IHS's admitting this
orbit space, most of which lie in non-K\"alerian closed symplectic
4-manifolds. Another simple example: Take 1/4 of the plane to be the orbit
space, so that the two half-lines correspond to elliptic singularities.
Equip this 1/4 plane with some affine structure, and we obtain an IHS lying
on a symplectic manifold diffeomorphic to ${\bf R}^4$. Most of the times
this ${\bf R}^4$ will be an exotic symplectic space.

To conclude this paper, we will give a physical example to illustrate our
results. Namely, we will discuss here one of the most famous IHS's ever
known: the Kovalevskaya's top. This top has attracted so many people, from
different points of view, since
its appearance in \cite{Kovalevskaya}. Kharlamov was first to study the
bifurcation of Liouville tori of this system (see e.g. \cite{Kharlamov}).
Here, for our purposes, we will follow Oshemkov \cite{Oshemkov}.
Remark that codimension 1 singularities of this system were reexamined by
Audin and Silhol \cite{AS} by the use of algebro-geometric methods.

Consider the Lie algebra $e(3)$ of the transformation group of the
three-dimensional Euclidean space. Let $S_1, S_2, S_3, R_1, R_2, R_3$ be a
system of coordinate functions on its dual $e(3)^*$, for which the
Lie-Poisson bracket is of the form:
$$ \{ S_i, S_j \} = \epsilon_{ijk} S_k, \{ R_i, R_j \} = 0,
\{ S_i, R_j \} = \{R_i, S_j\} = \epsilon_{ijk} R_k, $$
where $\{ i,j,k \} = \{1,2,3\}$ and \\
$\epsilon_{ijk} = \left\{ \begin{tabular}{l}
sign of the transposition of (i,j,k) if all the i,j,k are different \\
0 otherwise. \end{tabular} \right. $ \\

The above Lie-Poisson structure has 2 Casimir
functions (i.e. functions which commute with all other smooth
functions, with respect to the Lie-Poisson bracket):
$$f_1 = R_1^2 + R_2^2 + R^2_3, f_2 = S_1R_1 + S_2R_2 + S_3R_3 $$

The restriction of the Poisson structure on the
level sets $\{f_1 = const > 0, f_2 = const\}$ is non-degenerate. These level
sets are symplectic 4-manifolds, diffeomorphic to $TS^2$.

Motion of a rigid body with a fixed point in a gravity field
can often be written in the form
of a Hamiltonian system on $e(3)^*$, whence it admits $f_1, f_2$ as first
integrals, and is Hamiltonian when restricted to the above symplectic
4-manifolds. In case of the Kovalevskaya's top, the corresponding
Hamiltonian can be written in the form
$$ H = 1/2 (S_1^2 + S_2^2 + 2 S_3^2) + R_1$$

This Hamiltonian admits another first integral:
$$ K = (S_1^2/2 - S_2^2/2 - R_1)^2 + (S_1S_2 - R_2)^2 $$

Let us consider, for example, the restriction of this Kovalevskaya's system
to a particular symplectic 4-manifold $M^4 = \{ f_1 = 1, f_2 = g \}$, where
$g$ is a constant with $0 < |g| < 1 $. Then the bifurcation diagram of the
moment map $(H,K): M^4 \to {\bf R}^2$ is shown in the following figure.

\hspace{3.5cm} \epsfbox{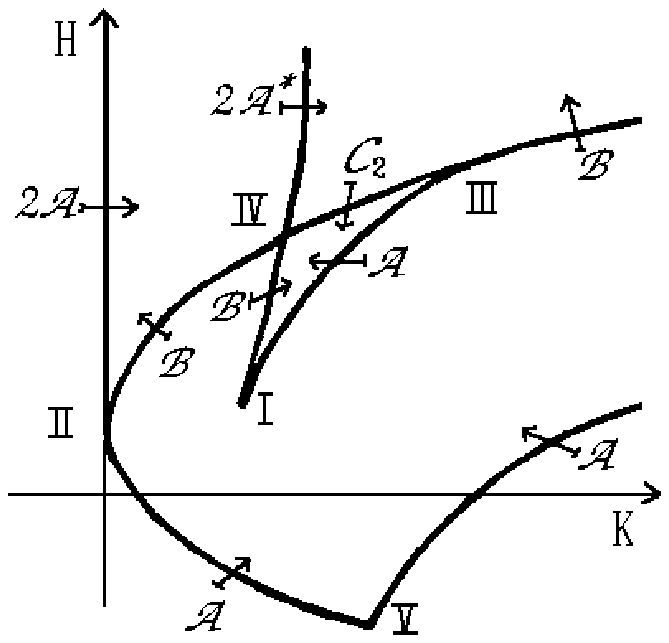} 

Together with
the bifurcation curves, in this figure we also show the type of codimension
1 nondegenerate singularities and point out simply-degenerate and
codimension 2 singularities. In the figure,
codimension 1 singularities are denoted by calligraphic letters, and these
notations and description are taken from \cite{BFM,Oshemkov}. In particular,
$\mathcal A$ means an elliptic codimension 1 singularity, $2 \mathcal A$ means a
disjoint union of 2 such singularities, $\mathcal B$ means the codimension 1
hyperbolic singularity for which the ${\bf T}^{n-1}$ action is free
(here $n-1 =1$), and which
contains only one hyperbolic orbit in the singular leaf. ${\mathcal A}^*$ is
obtained from $\mathcal B$ by a free ${\bf Z}_2$ action. ${\mathcal C}_2$
is the codimension 1 hyperbolic singularity for which the ${\bf T}^{n-1}$
action is free, and whose reduction is
a surface singularity which can be embedded in the 2-sphere in such a way
that the singular leaf is the union of two big circles (intersecting at 2
points).

The Kovalevskaya's top is nondegenerate but not strongly nondegenerate,
and there are 3 simply-degenerate singularities, which are denoted
in the above figure by I, II, III. Only singularity I (which corresponds to
the cusp in the bifurcation diagram) is simply-degenerate in the sense of
Lerman and Umanskii \cite{LU1}. Singularity V is codimension 2 elliptic, and
singularity IV is codimension 2 hyperbolic. Based on the structure of
codimension 1 singularities around singularity IV, and the decomposition
Theorem \ref{thm:adp}, we can easily compute singularity IV.
The result is: IV = $({\mathcal B} \times {\mathcal C}_2) / {\bf Z}_2$.



\end{document}